\documentclass[11pt]{article}

\usepackage{amsmath}
\usepackage{amsfonts}
\usepackage{amssymb}

\newcommand{\bB}{ {\mathbb B} }
\newcommand{\bC}{ {\mathbb C} }
\newcommand{\cD}{ {\cal D} }
\newcommand{\cE}{ {\cal E} }

\newcommand{\cH}{ {\cal H}}

\newcommand{\cK}{ {\cal K}}
\newcommand{\cM}{ {\cal M} }
\newcommand{\fM}{ {\mathfrak M} }
\newcommand{\bN}{ {\mathbb N} }

\newcommand{\fN}{ {\mathfrak N} }
\newcommand{\bR}{ {\mathbb R} }
\newcommand{\cR}{ {\cal R} }
\newcommand{\cS}{ {\cal S} }
\newcommand{\cT}{ {\cal T} }
\newcommand{\cV}{ {\cal V} }

\newcommand{\bZ}{ {\mathbb Z} }

\newcommand{\leqleq}{ \ll }
\newcommand{\geqgeq}{ \gg }
\newcommand{\fatone}{1} 
\newcommand{\UU}{\sqcup} 
\newcommand{\emb}{\Subset} 
\newcommand{\Parent}{ \mbox{Parent} } 
\newcommand{\piodd}{ \pi^{\mathrm{(odd)}} } 
\newcommand{\kpeven}{ K( \pi )^{\mathrm{(even)}} } 
\newcommand{\rhodd}{ \rho^{\mathrm{(odd)}} } 
\newcommand{\roeven}{ \rho^{\mathrm{(even)}} } 
\newcommand{\kroeven}{ K( \rho )^{\mathrm{(even)}} } 
\newcommand{\sigmaodd}{ \sigma^{\mathrm{(odd)}} } 
\newcommand{\taueven}{ \tau^{\mathrm{(even)}} } 

\newcommand{\Int}{ \mbox{Int} } 
\newcommand{\cf}{ \mbox{Cf} } 
\newcommand{\Reta}{ \mbox{Reta} } 
\newcommand{\card}{ \mbox{card} } 

\newcommand{\Dinfdiv}{ \cD_c^{\textrm{inf-div}} } 
\newcommand{\Rinfdiv}{ \cR_c^{\textrm{inf-div}} } 

\newcommand{\Dalg}{ \cD_{\mathrm{alg}} }
\newcommand{\ncpolk}{ \bC \langle X_1, \ldots , X_k \rangle }
\newcommand{\ncserk}{ \bC_0 \langle \langle z_1, \ldots ,z_k \rangle \rangle }

\newcommand{\toN}{ \lim_{N \to \infty} }

\newcommand{\freeplus}{ \boxplus }
\newcommand{\freetimes}{ \boxtimes }
\newcommand{\freestar}{ \framebox[7pt]{$\star$} }
\newcommand{\Uplus}{ \uplus }

\begin{document}

\title{\bf \boldmath{$\eta$}-series and a Boolean Bercovici--Pata \\
bijection for bounded \boldmath{$k$}-tuples }
\author{
Serban T. Belinschi 
\and Alexandru Nica \thanks{Research supported by a Discovery Grant of 
NSERC, Canada and by a PREA award from the province of Ontario.}  }

\date{ }

\maketitle

\vspace{-.5cm}

\begin{center}
Department of Pure Mathematics  

University of Waterloo 
\end{center}

\vspace{.5cm}

\begin{abstract}
Let $\cD_c (k)$ be the space of (non-commutative) distributions of 
$k$-tuples of selfadjoint elements in a $C^*$-probability space. On 
$\cD_c (k)$ one has an operation $\freeplus$ of free additive 
convolution, and one can consider the subspace $\Dinfdiv (k)$ of 
distributions which are infinitely divisible with respect to this 
operation. The linearizing transform for $\freeplus$ is the 
$R$-transform (one has $R_{\mu \freeplus \nu} = R_{\mu} + R_{\nu}$, 
$\forall \ \mu , \nu \in \cD_c (k)$). We prove that the set of 
$R$-transforms $\{ R_{\mu} \mid \mu \in \Dinfdiv (k) \}$ can also
be described as $\{ \eta_{\mu} \mid \mu \in \cD_c (k) \}$, where 
for $\mu \in \cD_c (k)$ we denote $\eta_{\mu} = M_{\mu}/ (1+M_{\mu})$, 
with $M_{\mu}$ the moment series of $\mu$. (The series $\eta_{\mu}$ 
is the counterpart of $R_{\mu}$ in the theory of Boolean convolution.) 
As a consequence, one can define a bijection 
$\bB : \cD_c (k) \to \Dinfdiv (k)$ via the formula 
\[
{ }  \hspace{3cm}
R_{\bB (\mu)} = \eta_{\mu}, \ \ \forall \, \mu \in \cD_c (k).  
\hspace{3.2cm}  \mbox{(I)}
\]
We show that $\bB$ is a multi-variable analogue of a bijection 
studied by Bercovici and Pata for $k=1$, and we prove a theorem 
about convergence in moments which parallels the Bercovici-Pata 
result. On the other hand we prove the formula
\[
{  }  \hspace{3cm}
\bB (\mu \freetimes \nu) = \bB ( \mu ) \freetimes \bB ( \nu ),
\hspace{3.5cm} \mbox{(II)}
\]
with $\mu, \nu$ considered in a space 
$\Dalg (k) \supseteq \cD_c (k)$ where the operation of free 
multiplicative convolution $\freetimes$ always makes
sense. An equivalent reformulation of (II) is that
\[
{  }  \hspace{2.5cm}
\eta_{\mu \freetimes \nu} = \eta_{\mu} \ \freestar \ \eta_{\nu},
\ \ \forall \, \mu , \nu \in \Dalg (k),
\hspace{2.4cm}   \mbox{(III)}
\]
where $\freestar$ is an operation on series previously studied by 
Nica and Speicher, and which describes the multiplication of free 
$k$-tuples in terms of their $R$-transforms. Formula (III) shows 
that, in a certain sense, $\eta$-series behave in the same way as 
$R$-transforms in connection to the operation of multiplication of 
free $k$-tuples of non-commutative random variables.
\end{abstract}

\newpage

\setcounter{section}{1}
\begin{center}
{\large\bf 1. Introduction}
\end{center}

The extent to which developments in free probability parallel phenomena 
from classical probability has exceeded by far what was originally expected 
in this direction of research. In particular there exists a well-developed 
theory of infinitely divisible distributions in the free sense; a few of the 
papers building this theory are \cite{V86}, \cite{BV92}, \cite{BP99} --
see also the section 2.11 of the survey \cite{V00} for more details. The 
Boolean Bercovici--Pata bijection is one of the results in this theory 
(cf \cite{BP99}, Section 6); it is a special bijection between the set 
of probability distributions on $\bR$ which are infinitely divisible with 
respect to free additive convolution (on one hand), and the set of all 
probability measures on $\bR$ (on the other hand).

In this paper we extend the Boolean Bercovici--Pata bijection to the 
multi-variable framework, in a context where we deal with $k$-tuples of
bounded random variables (or in other words, we deal with 
non-commutative multi-variable analogues for distributions with compact
support). The framework we consider is thus:
\[
\Bigl( \, \cD_c (k), \ \freeplus \, \Bigr),
\]
where $k$ is a positive integer, $\cD_c (k)$ is the set of linear 
functionals $\mu : \bC \langle X_1, \ldots , X_k \rangle \to \bC$ which 
appear as joint distribution for a $k$-tuple of selfadjoint elements in a
$C^*$-probability space, and $\freeplus$ is the operation of free additive 
convolution. This operation is considered in connection with the notion of 
free independence, and it encodes the distribution of the sum of two freely
independent $k$-tuples, in terms of the distributions of the two $k$-tuples. 
A more detailed review of $( \, \cD_c (k), \ \freeplus \, )$ and of the 
notations we are using in connection to it appears in the Section 4 of 
the paper. For a general introduction to the ideas of free probability we
refer to \cite{VDN92}. 

A distribution $\mu \in \cD_c (k)$ is said to be 
infinitely divisible with respect to $\freeplus$ if for every $N \geq 1$
there exists a distribution $\mu_N \in \cD_c (k)$ such that the 
$N$-fold $\freeplus$-convolution
$\mu_N \freeplus \mu_N \freeplus \cdots \freeplus \mu_N$
is equal to $\mu$. The set of distributions $\mu \in \cD_c (k)$ which 
have this property will be denoted by $\Dinfdiv (k)$.

The linearizing transform for $\freeplus$ is the {\em $R$-transform}.
The R-transform of a distribution $\mu \in \cD_c (k)$ is a power series 
$R_{\mu} \in \ncserk$, where $\ncserk$ denotes the set of power series with 
complex coefficients, and with vanishing constant term, in $k$ non-commuting 
indeterminates $z_1, \ldots , z_k$. The above mentioned linearization 
property is that
\begin{equation}   \label{eqn:1.1}
R_{\mu \freeplus \nu} = R_{\mu} + R_{\nu}, \ \ \forall \, 
\mu ,  \nu \in \cD_c (k).
\end{equation}
We denote by $\cR_c (k)$ 
the set of power series $f \in  \ncserk$ which appear as $R_{\mu}$ 
for some $\mu \in \cD_c (k)$; and, similarly, we use the notation 
$\Rinfdiv (k)$ for the set of series which appear as $R_{\mu}$ for a
distribution $\mu \in \Dinfdiv$. A distribution $\mu$ is always 
uniquely determined by its $R$-transform, so we are dealing in fact 
with two bijections,
\begin{equation}   \label{eqn:1.2}
\cD_c (k) \ni \mu \mapsto R_{\mu} \in \cR_c (k), \hspace{.5cm}
\Dinfdiv (k) \ni \mu \mapsto R_{\mu} \in \Rinfdiv (k).
\end{equation}

In this paper we put into evidence a commutative diagram where the vertical
arrows are the two bijections from (\ref{eqn:1.2}), and where the top 
horizontal arrow $\bB$ is a multi-variable counterpart of the 
Boolean Bercovici--Pata bijection:

\begin{equation}  \label{eqn:1.3}
\begin{matrix}
\cD_c (k) & \stackrel{\bB\ \ }{{\mbox{\huge{$\longrightarrow{}$}}}} & 
\Dinfdiv (k)\cr
& & & \cr
{{\mbox{\it\small R}}{\mbox{\huge{$\downarrow$}}}} &    \stackrel{\ \ \ 
\eta}{{\mbox{\huge$\searrow$}}}   & 
{{\mbox{\huge{$\downarrow$}}}{\mbox{\it\small R}}} \cr
   & & & \cr
\cR_c (k)  & \stackrel{\mbox{\scriptsize\Reta}\ \ 
}{{\mbox{\huge{$\longrightarrow{}$}}}} & \Rinfdiv (k)
\end{matrix}
\end{equation}

We observe the somewhat surprising occurrence in this diagram of another 
kind of transform, the {\em $\eta$-series}. For $\mu \in \cD_c (k)$, the 
$\eta$-series associated to $\mu$ is $\eta_{\mu}$ =
$M_{\mu} / (1+ M_{\mu}) \in \ncserk$, where $M_{\mu}$ denotes the moment 
series of $\mu$. The $\eta$-series is the linearizing transform for another 
kind of convolution on $\cD_c (k)$, the {\em Boolean convolution} $\uplus$.
The operation $\uplus$ is the counterpart of $\boxplus$ in connection to 
Boolean independence -- it encodes the distribution of the sum of two 
Boolean independent $k$-tuples, in terms of the distributions of the two 
$k$-tuples. The counterpart of Equation (\ref{eqn:1.1}) in the framework 
of $\uplus$ says that
\begin{equation}   \label{eqn:1.4}
\eta_{\mu \uplus \nu} = \eta_{\mu} + \eta_{\nu}, \ \ \forall \, 
\mu ,  \nu \in \cD_c (k).
\end{equation}
A distribution $\mu$ is uniquely determined by $\eta_{\mu}$, so if 
we denote $\cE_c (k) := \{ f \in \ncserk \mid $
$\exists \, \mu \in \cD_c (k)$ such that $\eta_{\mu} = f \}$, then 
we have a bijection
\begin{equation}   \label{eqn:1.5}
\cD_c (k) \ni \mu \mapsto \eta_{\mu}  \in \cE_c (k).
\end{equation}
Note that we have drawn this bijection along the diagonal of the 
diagram (\ref{eqn:1.3}); the justification for why we are allowed to do
so is given by the first part of the next theorem.

$\ $

{\bf Theorem 1.} {\em Let $k$ be a positive integer.

$1^o$ We have that $\Rinfdiv (k) \ = \ \cE_c (k)$
(equality of subsets of $\ncserk$). The target set of the bijection 
$\eta$ from (\ref{eqn:1.5}) can thus be regarded as $\Rinfdiv (k)$.

$2^o$ There exists a bijection
$\mbox{Reta} : \cR_c (k) \to \Rinfdiv (k)$ defined by the formula
\begin{equation}  \label{eqn:1.7}
\mbox{Reta} ( R_{\mu} ) = \eta_{\mu}, \ \ \forall \, \mu \in \cD_c (k).
\end{equation}
One has a purely combinatorial way of describing this bijection.
More precisely, there exists an explicit summation formula which 
gives the coefficients of Reta$(f)$ in terms of the coefficients of
$f$, where $f$ is an arbitrary series in $\cR_c (k)$. The summation 
formula is:
\begin{equation}  \label{eqn:1.8}
\cf_{(i_1, \ldots , i_n)} ( \Reta (f)) = 
\sum_{ \begin{array}{c} 
{\scriptstyle \pi \in NC(n),}  \\ 
{\scriptstyle \pi \leqleq 1_n} 
\end{array} } \ \cf_{(i_1, \ldots , i_n); \pi } (f) ,
\end{equation}
$\forall \, n \geq 1,$ $\forall \, 1 \leq i_1, \ldots , i_n \leq k$.
(The notations used in (\ref{eqn:1.8}) for coefficients of 
power series are detailed in Definition 3.2 below. The partial order 
``$\leqleq$'' on the set $NC(n)$ of non-crossing partitions of 
$\{ 1, \ldots , n \}$ is discussed in Section 2 below, starting with 
Definition 2.5.)

$3^o$ There exists a bijection $\bB : \cD_c (k) \to \Dinfdiv (k)$ 
which is determined by the formula
\begin{equation}  \label{eqn:1.9}
R_{\bB ( \mu )} = \eta_{\mu}, \ \ \forall \, \mu \in \cD_c (k).
\end{equation}
In the case $k=1$, $\bB$ coincides with the restriction of 
the Boolean Bercovici--Pata bijection to the set of compactly supported 
probability distributions on the real line. }

$\ $

When looking at the diagram (\ref{eqn:1.3}), one can say that the map 
Reta is the ``$R$-transform incarnation'' of the Boolean 
Bercovici--Pata 
bijection. It is remarkable that one can also describe Reta by the very 
explicit formulas (\ref{eqn:1.8}) and (\ref{eqn:1.7}) (where 
(\ref{eqn:1.7}) is the one which suggested the name ``Reta'' -- the 
transformation which ``converts $R$ to $\eta$''). 

We can supplement Theorem 1 with the following result, which is a
$k$-dimensional version of the equivalence $(2) \Leftrightarrow (3)$
in Theorem 6.1 of \cite{BP99} (and thus provides a more in-depth 
explanation for why $\bB$ of Theorem 1 is indeed a $k$-dimensional 
version of the corresponding bijection from \cite{BP99}).

$\ $

{\bf Theorem 1'.} {\em Let $k$ be a positive integer. Let 
$( \mu_N )_{N=1}^{\infty}$ be a sequence of distributions in $\cD_c (k)$,
and let $p_1< p_2< \cdots < p_N < \cdots$ be a sequence of positive 
integers. Then the following two statements are equivalent:

\vspace{6pt}

\begin{center}
(1) $\exists$ 
$\displaystyle\lim_{N \to \infty} \underbrace{\mu_N \freeplus \cdots 
\freeplus \mu_N}_{p_N \ times} \ =: \ \nu \in \Dinfdiv (k)$,

\vspace{10pt}

(2) $\exists$
$\displaystyle\lim_{N \to \infty} \underbrace{\mu_N \Uplus \cdots 
\Uplus \mu_N}_{p_N \ times} \ =: \ \mu \in \cD_c (k)$,
\end{center}

\vspace{6pt}

\noindent
where the limits  in (1), (2) are considered with respect to 
convergence in moments. Moreover, if (1) and (2) 
hold, then the resulting limits $\mu , \nu$ are connected by the 
formula $\bB ( \mu ) = \nu$, where $\bB$ is the bijection from
Theorem 1. }

$\ $

We next proceed to presenting the second main result of this paper, 
which concerns a surprising property
of the Boolean Bercovici--Pata bijection, in connection to the operation
of free {\em multiplicative} convolution. This result takes place in a 
purely algebraic framework, and in order to present it we will move from 
$\cD_c (k)$ to the larger set $\Dalg (k)$ of distributions of $k$-tuples
in arbitrary (purely algebraic) non-commutative probability spaces. 
$\Dalg (k)$ consists in fact of {\em all} linear functionals 
$\mu : \ncpolk \to \bC$ which satisfy the normalization condition 
$\mu ( \fatone ) = 1$.

The commutative diagram (\ref{eqn:1.3}) has a ``simplified'' version
living in the algebraic framework of $\Dalg (k)$. Indeed, for 
$\mu \in \Dalg (k)$ it is still possible to define the $R$-transform 
$R_{\mu}$, and every series $f \in \ncserk$ can be written uniquely as 
$R_{\mu}$ for some $\mu \in \Dalg (k)$. Moreover, every $\mu \in \Dalg (k)$
is (trivially) infinitely divisible in this purely algebraic framework;
hence the two bijections displayed in (\ref{eqn:1.2}) are now both 
replaced by the bijection
\begin{equation}  \label{eqn:1.10}
\Dalg (k) \ni \mu \mapsto R_{\mu} \in \ncserk .
\end{equation}
On the other hand, for $\mu \in \Dalg (k)$ one can define the $\eta$-series 
$\eta_{\mu}$, and every $f \in \bC \langle \langle z_1, \ldots ,$
$z_k \rangle \rangle$ can be written uniquely as 
$\eta_{\mu}$ for some $\mu \in \Dalg (k)$; so we also have a bijection
\begin{equation}  \label{eqn:1.105}
\Dalg (k) \ni \mu \mapsto \eta_{\mu} \in \ncserk ,
\end{equation} 
which is the counterpart of the bijection from (\ref{eqn:1.5}). This 
leads to the diagram
\begin{equation}  \label{eqn:1.11}
\begin{matrix}
\Dalg (k) & \stackrel{\bB\ \ }{{\mbox{\huge{$\longrightarrow{}$}}}} & 
\Dalg (k) \cr
& & & \cr
{{\mbox{\it\small R}}{\mbox{\huge{$\downarrow$}}}} &    \stackrel{\ \ \ 
\eta}{{\mbox{\huge$\searrow$}}}   & 
{{\mbox{\huge{$\downarrow$}}}{\mbox{\it\small R}}} \cr
   & & & \cr
\ncserk   & \stackrel{\mbox{\scriptsize\Reta}\ \ 
}{{\mbox{\huge{$\longrightarrow{}$}}}} & \ncserk
\end{matrix}
\end{equation}
where the vertical arrows are from (\ref{eqn:1.10}), the diagonal is from 
(\ref{eqn:1.105}), and the horizontal arrows $\bB$ and $\Reta$ are defined 
via the requirement that the diagram is commutative.

On the space $\Dalg (k)$ we can define an operation of 
{\em free multiplicative convolution} $\freetimes$, as follows. Given 
$\mu, \nu \in \Dalg (k)$, one can always find random variables 
$x_1, \ldots , x_k, y_1, \ldots , y_k$ in a non-commutative probability 
space $( \cM, \varphi )$ such that the joint distribution of the $k$-tuple 
$x_1, \ldots , x_k$ is equal to $\mu$, the joint distribution of the 
$k$-tuple $y_1, \ldots , y_k$ is equal to $\nu$, and such that 
$\{ x_1, \ldots , x_k \}$ is freely independent from 
$\{ y_1, \ldots , y_k \}$ in $( \cM , \varphi )$. The joint distribution 
of the $k$-tuple $x_1y_1, \ldots , x_ky_k$ turns out to depend only 
on $\mu$ and $\nu$; and the free multiplicative convolution 
$\mu \freetimes \nu$ is equal, by definition, to the joint 
distribution of $x_1y_1, \ldots , x_ky_k$. (Note:
what makes this definition not to work in the framework of $\cD_c (k)$
is that, even if we assume that all of $x_1, \ldots , x_k,y_1, \ldots , y_k$
are selfadjoint elements in a $C^*$-probability space, the products 
$x_1y_1, \ldots , x_ky_k$ will no longer be selfadjoint, in general.)

By using this terminology, our second theorem is then stated as follows.

$\ $

{\bf Theorem 2.} 
{\em The bijection $\bB$ from the commutative diagram (\ref{eqn:1.11}) 
is a homomorphism for $\freetimes$.
That is, we have 
\begin{equation}  \label{eqn:1.12}
\bB ( \mu \freetimes \nu ) = \bB ( \mu ) \freetimes \bB ( \nu ), \ \ 
\forall \, \mu , \nu \in \Dalg (k).
\end{equation}  }

$\ $

It is also worth recording how Theorem 2 looks like when it is re-phrased 
in terms of $R$-transforms. This re-phrasing involves an operation 
$\freestar$, called {\em boxed convolution}, on the space of series 
$\ncserk$. One way of defining $\freestar$ is via the equation
\begin{equation}  \label{eqn:1.13}
R_{\mu \freetimes \nu} = R_{\mu} \ \freestar \ R_{\nu},
\ \ \, \forall \, \mu , \nu \in \Dalg (k).
\end{equation}
This equation says that $\freestar$ is the ``incarnation of $\freetimes$''
obtained when one moves from $\Dalg (k)$ to $\ncserk$ via the bijection 
(\ref{eqn:1.10}). (Or, if we recall how $\freetimes$ is defined, we can
say that the job of $\freestar$ is to describe the multiplication of 
freely independent $k$-tuples, in terms of their $R$-transforms.) On the
other hand, the operation $\freestar$ can also be introduced in a purely 
combinatorial way -- one has explicit formulas giving the coefficients 
of $f \ \freestar \ g$ in terms of the coefficients of $f$ and of $g$, via 
summations over non-crossing partitions. The explicit formulas for the
coefficients of $f \ \freestar \ g$ will be reviewed in the Section 7 of the 
paper; for more details on $\freestar$ (including the explanation of why 
it is justified to use the name ``convolution'' for this operation) we refer
to \cite{NS06}, Lectures 17 and 18.

The reformulation of Theorem 2 in terms of transforms goes as follows.

$\ $

{\bf Theorem 2'.} {\em The $\eta$-series satisfies the relation
\begin{equation}  \label{eqn:1.14}
\eta_{\mu \freetimes \nu} = \eta_{\mu} \ \freestar \ \eta_{\nu},
\ \ \forall \, \mu , \nu \in \Dalg (k).
\end{equation}  }

$\ $

The equivalence of Theorems 2 and 2' is immediate. For instance if 
we assume Theorem 2', then Theorem 2 is obtained as follows: for 
every $\mu , \nu \in \Dalg (k)$ we have that
\begin{align*}
R_{\bB ( \mu \freetimes \nu )} 
& = \eta_{\mu \freetimes \nu} \ \mbox{(from diagram (\ref{eqn:1.11}))} \\
& = \eta_{\mu} \ \freestar \ \eta_{\nu} \ \mbox{(by Theorem 2')}       \\
& = R_{\bB ( \mu )} \ \freestar \ R_{\bB ( \nu )} 
                              \ \mbox{(from diagram (\ref{eqn:1.11}))} \\
& = R_{\bB ( \mu ) \freetimes \bB ( \nu )}  \
                                \mbox{(by Equation (\ref{eqn:1.13})).}
\end{align*}
So $\bB ( \mu \freetimes \nu )$ and $\bB ( \mu ) \freetimes \bB ( \nu )$
have the same R-transform, and these two distributions must therefore
be equal to each other.

It is interesting to compare the Equation (\ref{eqn:1.14}) in Theorem 2' 
with the quite similarly looking Equation (\ref{eqn:1.13}) which precedes 
the theorem.
We see here that the operation of boxed convolution $\freestar$ also 
pops up as the ``incarnation of $\freetimes$'' when one moves from 
$\Dalg (k)$ to $\ncserk$ via the bijection (\ref{eqn:1.105}), 
$\mu \mapsto \eta_{\mu}$, 
in lieu of the bijection $\mu \mapsto R_{\mu}$ from (\ref{eqn:1.10}).
(The bijection in (\ref{eqn:1.105}) is quite a bit easier to work with 
than the one in (\ref{eqn:1.10}) -- see the discussion in Proposition 3.5
and Remark 3.6 below.)

Let us also mention here that in the case when $k=1$, one of the usual 
ways of looking at $\freetimes$ is by viewing it as an operation on the set 
of probability distributions with support (not necessarily compact) 
contained in $[0, \infty )$. In this framework, a further discussion around 
the $\freetimes$-multiplicativity of the bijection $\bB$ is made in 
\cite{BN06} (by using complex analysis methods specific to the case $k=1$,
which also cover the situation of unbounded supports).

$\ $

We conclude this introductory section by describing how the rest of
the paper is organized. In the above discussion it was more relevant 
to consider first the framework of $\cD_c (k)$, but for a more detailed 
presentation it is actually better to first clarify the simpler algebraic
framework of $\Dalg (k)$. This is done in Section 3 of the paper, following
to a review of some basic combinatorial structures done in Section 2. 
In Section 4 we give a more 
detailed introduction to $\cD_c (k)$ and to the maps involved in 
the commutative diagram (\ref{eqn:1.3}), and then in Section 5 we give
the proofs of Theorems 1 and 1'. In Section 6 we return to the algebraic 
framework and present the result from the combinatorics of non-crossing 
partitions (Corollary 6.11) which lies at the core of our Theorems 2 
and 2'. Finally, Section 7 is devoted to presenting the proofs of 
Theorems 2 and 2'.

$\ $

$\ $

\setcounter{section}{2}
\begin{center}
{\large\bf 2. Some basic combinatorial structures}
\end{center}

The first part of this section gives a very concise review (intended mostly 
for setting notations) of non-crossing partitions. For a more
\setcounter{equation}{0}
detailed introduction to these partitions, and on how they are used 
in free probability, we refer to \cite{NS06}, Lectures 9 and 10.

$\ $

{\bf 2.1 Remark} {\em (review of $NC(n)$).}

Let $n$ be a positive integer and let 
$\pi = \{ B_{1} , \ldots , B_{p} \}$ be a partition of $\{ 1, \ldots ,n \}$ 
-- i.e. $B_{1} , \ldots , B_{p}$ are pairwise disjoint non-void sets 
(called the {\bf blocks} of $\pi$), and $B_{1} \cup \cdots \cup B_{p}$ = 
$\{ 1, \ldots , n \}$. We say that $\pi$ is {\bf non-crossing} if for 
every $1 \leq i < j < k < l \leq n$ such that $i$ is in the same
block with $k$ and $j$ is in the same block with $l$, it necessarily
follows that all of $i,j,k,l$ are in the same block of $\pi$.
The set of all non-crossing partitions of $\{ 1, \ldots , n \}$ will be
denoted by $NC(n).$ On $NC(n)$ we consider the partial order given 
by {\bf reversed refinement}: for $\pi , \rho \in NC(n)$, we write
``$\pi \leq \rho$'' to mean that every block of $\rho$ is a union of
blocks of $\pi$. (In this paper
we will use more than one partial order on $NC(n)$, but ``$\leq$'' will 
be always reserved for reversed refinement order.)

For $\pi \in NC(n)$, the number of blocks of $\pi$ will be denoted
by $| \pi |$. The minimal and maximal element of $( NC(n), \leq )$ are 
denoted by $0_n$ (the partition of $\{ 1, \ldots , n \}$ into $n$ 
singletons) and respectively $1_n$ (the partition of $\{ 1, \ldots , n \}$ 
into one block).

A partition $\pi \in NC(n)$ has an {\bf associated permutation} of 
$\{ 1, \ldots , n \}$, which will be denoted by $P_{\pi}$. The permutation
$P_{\pi}$ is defined by the prescription that for every block 
$B = \{ b_1, \ldots , b_m \}$ of $\pi$, with $b_1 < \cdots < b_m$, one 
creates a cycle of $P_{\pi}$, as follows:
\[
P_{\pi} (b_1) = b_2, \ldots , P_{\pi} (b_{m-1}) = b_m, P_{\pi} (b_m) = b_1.
\]

$\ $

{\bf 2.2 Remark} {\em (review of the Kreweras complementation map).} 

This is a special order-reversing bijection $K: NC(n) \to NC(n)$. One way 
of describing how it works (which is actually the original definition from 
\cite{K72}) goes by using partitions of $\{ 1, \ldots , 2n \}$.

Let $\pi$ and $\rho$ be two partitions of $\{ 1,\ldots , n \}$. We will 
denote by
\[
\piodd \UU \roeven
\]
the partition of $\{ 1, \ldots , 2n \}$ which is obtained when one turns
$\pi$ into a partition of $\{ 1,3, \ldots ,$
$2n-1 \}$ and one turns $\rho$
into a partition of $\{ 2,4, \ldots , 2n \}$, in the canonical way. That is,
$\piodd \UU \roeven$ has blocks of the form $\{ 2a-1 \mid a \in A \}$ where 
$A$ is a block of $\pi$, and has blocks of the form $\{ 2b \mid b \in B \}$ 
where $B$ is a block of $\rho$. 

A partition $\theta$ of $\{ 1, \ldots , 2n \}$ is said to be 
{\bf parity-preserving} if every block of $\theta$ either is contained in
$\{ 1, 3 , \ldots , 2n-1 \}$ or is contained in $\{ 2,4, \ldots , 2n \}$.
The partitions of the form $\piodd \UU \roeven$ introduced 
above are parity-preserving; and conversely, every parity-preserving
partition $\theta$ of $\{ 1, \ldots , 2n \}$ is of the form 
$\piodd \UU \roeven$ for some uniquely determined partitions $\pi, \rho$
of $\{ 1, \ldots , n \}$.

The requirement that $\pi$ and $\rho$ are in $NC(n)$ is clearly necessary
but not sufficient in order for $\piodd \UU \roeven$ to be in $NC(2n)$.
If we fix $\pi \in NC(n)$ then the set
\[
\{ \rho \in NC(n) \mid \piodd \UU \roeven \in NC(2n) \}
\]
turns out to contain a largest partition $\rho_{\mathrm{max}}$, which is 
called the {\bf Kreweras complement} of $\pi$ and is denoted by $K( \pi )$.
So $K( \pi )$ is defined by the requirement that for $\rho \in NC(n)$ 
we have:
\begin{equation}  \label{eqn:2.01}
\piodd \UU \roeven \in NC(2n) \ \Leftrightarrow \ \rho \leq K( \pi ). 
\end{equation}

It is easily verified that $\pi \mapsto K( \pi )$ is indeed 
an order-reversing bijection from $NC(n)$ to itself. Another feature of 
Kreweras complementation which is worth recording is that 
\begin{equation}   \label{eqn:2.02}
|\pi| + |K( \pi )| = n+1, \ \ \forall \, \pi \in NC(n).
\end{equation}

$\ $

{\bf 2.3 Remark} {\em (Kreweras complementation via permutations).}

A convenient way of describing Kreweras complements is by using the 
permutations associated to non-crossing partitions. Indeed, the permutation
$P_{K( \pi )}$ associated to the Kreweras complement of $\pi$ turns out to
be given by the neat formula
\begin{equation}  \label{eqn:2.1}
P_{K( \pi )} = P_{\pi}^{-1} P_{1_n}, \ \ \ \pi \in NC(n).
\end{equation}
(Note that the permutation $P_{1_n}$ associated to the maximal partition
$1_n \in NC(n)$ is just the cycle 
$1 \mapsto 2 \mapsto \cdots \mapsto n \mapsto 1$.)

The formula (\ref{eqn:2.1}) can be extended in order to cover the concept
of {\bf relative Kreweras complement} of $\pi$ in $\rho$, for 
$\pi, \rho \in NC(n)$ such that $\pi \leq \rho$. This relative Kreweras 
complement is a partition in 
$NC(n)$, which will be denoted by $K_{\rho} ( \pi )$, and which is 
uniquely determined by the fact that the permutation associated to it is
\begin{equation}  \label{eqn:2.2}
P_{K_{\rho} ( \pi )} = P_{\pi}^{-1} P_{\rho}.
\end{equation}
Clearly, the Kreweras complementation map $K$ discussed above
is the relative complementation with respect to the maximal 
element $1_n$ of $NC(n)$.

It can be shown that, for a fixed $\rho \in NC(n)$, the map 
$\pi \mapsto K_{\rho} ( \pi )$
is an order-reversing bijection from 
$\{ \pi \in NC(n) \mid \pi \leq \rho \}$ onto itself.
It can also be shown that 
\begin{equation}  \label{eqn:2.3}
\pi \leq \rho_1  \leq \rho_2 \mbox{ in $NC(n)$ } \ \Rightarrow \ 
K_{\rho_1} ( \pi ) \leq K_{\rho_2} ( \pi ).
\end{equation}
For proofs of these facts, and for more details on relative 
Kreweras complements we refer to \cite{NS06}, Lecture 18.

$\ $

Besides $NC(n)$, we will also use the partially 
ordered set of {\em interval partitions}.

$\ $

{\bf 2.4 Remark} {\em (review of $Int (n)$).} 

A partition $\pi$ of $\{ 1, \ldots , n \}$ is 
said to be an {\bf interval partition} if every block $B$ of $\pi$ is of the 
form $B = [i,j] \cap \bZ$ for some $1 \leq i \leq j \leq n$. The set of all 
interval partitions of $\{ 1, \ldots , n \}$ will be denoted by $\Int (n)$.
It is clear that $\Int (n) \subseteq NC(n)$, but it is in fact customary
to view $( \Int (n), \leq )$ as a partially ordered set on its own (where 
``$\leq$'' still stands for the reversed refinement order on partitions).
The enumeration arguments related to $\Int (n)$ are often simplified by 
the fact that we have a natural bijection between $\Int (n)$ and the 
collection $2^{\{ 1, \ldots , n-1 \} }$ of all subsets of 
$\{ 1, \ldots , n-1 \}$; this bijection maps $\pi \in \Int (n)$ to the 
set
\[
\left\{ m \ \begin{array}{ll}
\vline  & 1 \leq m \leq n-1, \mbox{and there exists a}  \\
\vline  & \mbox{block $B$ of $\pi$ such that $\max (B) =m$}
\end{array}  \right\} . 
\] 
Moreover, this bijection is a poset isomorphism, if one endows 
$2^{\{ 1, \ldots , n-1 \} }$ with the partial order given by 
reversed inclusion.

$\ $

We now move to introduce another partial order on $NC(n)$;
this is not part of the usual lingo related to this topic, but will 
turn out to be essential for the developments shown in the present paper.

$\ $

{\bf 2.5 Definition.} Let $n$ be a positive integer, and let $\pi$ and
$\rho$ be two partitions in $NC(n)$. We will write
``$\pi \leqleq \rho$'' to mean that $\pi \leq \rho$ and that, in addition, 
the following condition is fulfilled:
\begin{equation}  \label{eqn:3.1}
\left\{  \begin{array}{l}
\mbox{For every block $C$ of $\rho$ there exists a block}  \\
\mbox{$B$ of $\pi$ such that $\min (C), \max (C) \in B$.}
\end{array}  \right.
\end{equation}

$\ $

{\bf 2.6 Remark.} 
$1^o$ Let $\pi , \rho \in NC(n)$ be such that 
$\pi \leqleq \rho$. Let $C$ be a block of $\rho$, and let $B$ be the 
block of $\pi$ which contains $\min (C)$ and $\max (C)$. Then 
$B \subseteq C$ (because $B$ has to be contained in a block of $\rho$,
and this block can only be $C$), and we must have 
\begin{equation}  \label{eqn:3.2}
\min (B) = \min (C), \ \ \max (B) = \max (C). 
\end{equation}

$2^o$ It is immediately verified that $\leqleq$ is indeed a 
partial order relation on $NC(n)$. It is much coarser than the reversed 
refinement order. For instance, the inequality $\pi \leqleq 1_n$ is not 
holding for all $\pi \in NC(n)$, but it rather amounts to the condition 
that the numbers $1$ and $n$ belong to the same block of $\pi$. At the 
other end of $NC(n)$, the inequality $\pi \geqgeq 0_n$ can only take place
when $\pi = 0_n$. (While looking at these trivial examples, let us also 
note that the partial order $\leqleq$ does not generally behave 
well under taking Kreweras complements.)

$3^o$ Let $\rho = \{ C_1, \ldots , C_p \}$ be a fixed partition in $NC(n)$. 
For every $1 \leq q \leq p$ such that $|C_q| \geq 3$, let us split the 
block $C_q$ into the doubleton $\{ \min (C_q), \max (C_q) \}$ and 
$| C_q | - 2$ singletons; when doing this for all $q$ we obtain a 
partition $\rho_0 \leq \rho$ in $NC(n)$, such that all the blocks of 
$\rho_0$ have either 1 or 2 elements. From Definition 2.5 it is clear that
for $\pi \in NC(n)$ we have:
\begin{equation}  \label{eqn:3.3}
\pi \leqleq \rho \ \Leftrightarrow \ \rho_0 \leq \pi \leq \rho .
\end{equation}
Consequently, the set $\{ \pi \in NC(n) \ | \ \pi \leqleq \rho \}$
is just the interval $[ \rho_0 , \rho ]$ (with respect to reversed 
refinement order) of $NC(n)$, and in order to describe
it one can use the nice structure of such intervals of $NC(n)$ -- as 
presented for instance in \cite{NS06}, Lecture 9.

$4^o$ Let $\pi$ be a fixed partition in $NC(n)$. In contrast to what 
was observed in the preceding part of this remark, the set 
$\{ \rho \in NC(n) \ | \ \rho  \geqgeq \pi \}$
isn't generally an interval with respect to reversed refinement order.
This set has nevertheless nice enumerative properties, which will be 
described in Proposition 2.13 below. In Proposition 2.13 we will 
use a few basic facts concerning the nested structure of the blocks 
of a non-crossing partition, and we start by presenting these facts.

$\ $

{\bf 2.7 Definition.} 
Let $n$ be a positive integer and let $A,B$ be 
two non-empty subsets of $\{ 1, \ldots , n \}$. If 
$\min (A) \leq \min (B)$ and $\max (A) \geq \max (B)$, then we will
say that $A$ {\bf embraces} $B$, and write $A \emb B$. 

$\ $

{\bf 2.8 Definition.} 
$1^o$ Let $\pi$ be a partition in $NC(n)$, and let $A$ be a block 
of $\pi$. If there is no block $B$ of $\pi$ such that
$\min (B) < \min (A) \leq \max (A) < \max (B)$,
then we say that $A$ is an {\bf outer block} of $\pi$.

$2^o$ For $\pi \in NC(n)$, the number of outer blocks of $\pi$ will be 
denoted as $| \pi |_{\mathrm{out}}$.

$\ $

{\bf 2.9 Remark.} Let $\pi$ be a partition in $NC(n)$.

$1^o$ It is clear that the set of blocks of $\pi$ is partially
ordered by embracing (where we stipulate that the block $A$ is 
smaller than the block $B$ with respect to this partial order 
if and only if $A \emb B$). It is also clear that the outer blocks
are precisely the minimal blocks of $\pi$ with respect to this 
partial order.

$2^o$ Concerning the outer blocks of $\pi$, it is immediate that:

(a) The block of $\pi$ which contains the number 1 is outer.

(b) If $B$ is an outer block of $\pi$ such that $\max (B) <n$, then 
there exists an outer block $B'$ of $\pi$ such that 
$\min (B') = \max (B)+1$.

Hence if we denote $| \pi |_{\mathrm{out}} =:r$ and if we list the 
outer blocks of $\pi$ as $B_1, B_2, \ldots , B_r$,
in increasing order of their minimal elements, then we have
\[
\min (B_1) =1, \ \min (B_2) = \max (B_1)+1, \ldots ,
\min (B_r) = \max (B_{r-1})+1, \ \max (B_r) =n.
\]

$3^o$ Let $B_1, \ldots , B_r$ be as in the preceding part of this
remark, and let $\rho$ be the interval partition with blocks 
$[ \min (B_i) , \max (B_i) ] \cap \bZ$, $1 \leq i \leq r$. It is 
immediate that $\rho \geq  \pi$ and that $\rho$ is the smallest 
interval partition (in the sense of reversed refinement order) 
which satisfies this inequality.

$\ $

{\bf 2.10 Proposition.} Let $\pi$ be a partition in $NC(n)$.

$1^o$ Let $A$ and $B$ be two distinct blocks of $\pi$. We have
\begin{equation}  \label{eqn:8.1}
A \emb B \ \Leftrightarrow \ \Bigl( \, \exists \, a_1,a_2 \in A
\mbox{ and } b \in B \mbox{ such that } a_1 < b < a_2 \, \Bigr) .
\end{equation}

$2^o$ Let $A_1, A_2,B$ be blocks of $\pi$ such that $A_1 \emb B$ 
and $A_2 \emb B$. Then either $A_1 \emb A_2$ or $A_2 \emb A_1$.

$\ $

{\bf Proof.} $1^o$ is immediate, and left as exercise. For $2^o$ 
we observe first that the intersection 
$[ \min (A_1) , \max (A_1) ] \cap [ \min (A_2) , \max (A_2) ]$ is 
non-empty, as it contains $B$. If it was not true that one of the 
intervals $[ \min (A_1), \max (A_1) ]$,
$[ \min (A_2) , \max (A_2) ]$ is contained in the other, we would 
obtain either that $\min (A_1) < \min (A_2) < \max (A_1) < \max (A_2)$
or that $\min (A_2) < \min (A_1) < \max (A_2) < \max (A_1)$, 
contradicting the fact that $\pi$ is non-crossing.  {\bf QED}

$\ $

{\bf 2.11 Remark.} Let $\pi$ be a partition in $NC(n)$. We consider 
the partial order given by embracing on the set of blocks of $\pi$, 
and we consider the so-called {\em Hasse diagram} for this partial
order. The Hasse diagram is, by definition, a graph which has vertex
set equal the set of blocks of $\pi$, and has an edge connecting 
two blocks $A_1 \neq A_2$ when one of them embraces the other
(say for instance that $A_1 \emb A_2$) and there is no third block 
of $\pi$ which lies strictly between them (no
$A \neq A_1, A_2$ such that $A_1 \emb A \emb A_2$).
It is instructive to note that, as an immediate consequence of 
Proposition 2.10.2, the Hasse diagram we just described is a 
{\em forest} -- that is, each of its connected components is a 
{\em tree} (a graph without circuits).

Let us also recall here the following concept: a forest is said to
be {\em rooted} if one special vertex (a ``root'') has been chosen
in each of its connected components. There exists a natural way of
rooting the above Hassse diagram because, as immediately seen, each
of its connected components contains precisely one outer block of 
$\pi$. Thus we can view the Hasse diagram as a rooted forest, where 
the outer blocks are the roots.

$\ $

We now return to the partial order $\leqleq$ on $NC(n)$ that was 
introduced in Definition 2.5.

$\ $

{\bf 2.12 Remark.} 
Let $\pi , \rho$ be two partitions in $NC(n)$ such that
$\pi \leqleq \rho$. Let $C$ be a block of $\rho$ and let $B$ be 
the unique block of $\pi$ which has $\min (B) = \min (C)$ and 
$\max (B) = \max (C)$ (cf Remark 2.6.1). It is an easy exercise, 
left to the reader, to check that 
\begin{equation}  \label{eqn:2.110}
\left(  \begin{array}{c}
\mbox{$C$ is an outer}   \\
\mbox{block of $\rho$}  
\end{array}  \right) \ \Leftrightarrow \
\left(  \begin{array}{c}
\mbox{$B$ is an outer}   \\
\mbox{block of $\pi$}  
\end{array}  \right)  .
\end{equation}
Moreover, it is easily seen that in this way we obtain a canonical 
one-to-one correspondence between the outer blocks of $\rho$ and 
the outer blocks of $\pi$.

As a consequence of the above, we get that for a given 
$\pi \in NC(n)$, the number of blocks of any partition 
$\rho \in NC(n)$ such that $\rho \geqgeq \pi$ is bounded by 
the inequalities
\begin{equation}
| \pi |_{\mathrm{out}} \leq | \rho | \leq | \pi |.
\end{equation}
(The first of the two inequalities holds because
$| \pi |_{\mathrm{out}} = | \rho |_{\mathrm{out}} \leq | \rho |$,
while the second one follows directly from the fact that 
$\rho \geq \pi$.)

$\ $

{\bf 2.13 Proposition.} Let $\pi$ be a partition in $NC(n)$. For 
every integer $p$ satisfying 
$| \pi |_{\mathrm{out}} \leq p \leq | \pi |$, we have that:
\begin{equation} \label{eqn:2.12}
\card \Bigl\{ \rho \in NC(n) \mid \rho \geqgeq \pi \mbox{ and } 
| \rho | =p \Bigr\} 
= \left(  \begin{array}{c}
| \pi | - | \pi |_{\mathrm{out}}  \\
   p    - | \pi |_{\mathrm{out}} 
\end{array}  \right) .
\end{equation}

$\ $

{\bf Proof.} Consider the following two conditions which a set 
$\fM$ of blocks of $\pi$ may fulfil.
\[
\mbox{(C)} \ \  \left\{ \begin{array}{l} 
\mbox{{\em ``Convexity condition''.} Whenever $A_1, A_2, A_3$ 
                         are blocks of $\pi$ such}                \\
\mbox{ that $A_1 \emb A_2 \emb A_3$ and such that 
                        $A_1, A_3 \in \fM$, it follows}           \\
\mbox{ that $A_2 \in \fM$ as well.}
\end{array}  \right.
\]

\[
\mbox{(Min)} \ \  \left\{ \begin{array}{l} 
\mbox{{\em ``Minimal element condition''.} There exists a 
                    (necessarily unique)}                     \\
\mbox{ block $B \in \fM$ such that $B \emb A$, 
                       $\forall \, A \in \fM$.}              
\end{array}  \right.
\]
Throughout this proof we will use the ad-hoc term of 
{\em ``CMin set of blocks''} for a set $\fM$ of blocks of $\pi$ 
which fulfills both the conditions (C) and (Min). Moreover, suppose 
that $\fM_1, \ldots , \fM_p$ are CMin sets of blocks of $\pi$ such 
that $\fM_i \cap \fM_j = \emptyset$ for $i \neq j$, and such that 
$\fM_1 \cup \cdots \cup \fM_p$ contains all blocks of $\pi$;
then we will refer to $\{ \fM_1, \ldots , \fM_p \}$ by calling it 
a {\em ``CMin decomposition for the set of blocks of $\pi$''}.

The relevance of CMin sets of blocks in this proof comes from the 
following fact.

\vspace{10pt}

{\em Fact 1. Let $\rho$ be a partition in $NC(n)$ such that 
$\rho \geqgeq \pi$, let $C$ be a block of $\rho$, and denote
$\fM := \{ A \mid A$ block of $\pi$, $A \subseteq C \}$.
Then $\fM$ is a CMin set of blocks of $\pi$. }

\vspace{10pt}

The verification of Fact 1 is immediate. For instance in order to 
verify that $\fM$ satisfies the condition (C), one proceeds as 
follows. Let $A_1, A_2, A_3$ be blocks of $\pi$ such 
that $A_1 \emb A_2 \emb A_3$ and such that $A_1, A_3 \in \fM$. 
Assume by contradiction that $A_2 \not\in \fM$, hence that 
$A_2 \subseteq C'$ for some block $C'$ of $\rho$ where 
$C' \neq C$. From the hypothesis $A_1 \emb A_2 \emb A_3$ we deduce
that $\min (A_1) < \min (A_2) < \min (A_3) < \max (A_2)$;
since $\min (A_1) , \min (A_3) \in C$ and 
$\min (A_2) , \max (A_2) \in C'$, we have thus obtained a crossing
between the blocks $C$ and $C'$ of $\rho$ -- contradiction. 

An immediate consequence of Fact 1 is that we have:

\vspace{10pt}

{\em Fact 2. Let $\rho = \{ C_1, \ldots , C_p \}$ be a partition 
in $NC(n)$ such that $\rho \geqgeq \pi$, and for every 
$1 \leq i \leq p$ let us denote
$\fM_i := \{ A \mid \mbox{$A$ block of $\pi$, $A \subseteq C_i$} \}$.
Then $\{ \fM_1, \ldots , \fM_p \}$ is a CMin decomposition of the 
set of blocks of $\pi$. }

\vspace{10pt}

On the other hand we have a converse of the Fact 2, stated as 
follows.

\vspace{10pt}

{\em Fact 3. Let $\{ \fM_1, \ldots , \fM_p \}$ be a CMin 
decomposition of the set of blocks of $\pi$, and for every
$1 \leq i \leq p$ let us denote
\begin{equation}  \label{eqn:2.13}
C_i \ := \ \cup_{A \in \fM_i} \ A.
\end{equation}
Consider the partition $\rho = \{ C_1, \ldots , C_p \}$ 
of $\{ 1, \ldots , n \}$. Then $\rho \in NC(n)$, and 
$\rho \geqgeq \pi$. }

\vspace{10pt}

{\em Verification of Fact 3.} The only non-trivial point in the 
statement of Fact 3 is that the partition $\rho$ is non-crossing.
In order to verify this, let us fix 
$1 \leq a_1 < a_2 < a_3 < a_4 \leq n$ such that $a_1, a_3 \in C_i$
and $a_2, a_4 \in C_j$ for some $1 \leq i,j \leq p$. We want to 
prove that $i=j$.

Let $A'$ and $A''$ denote the blocks of $\pi$ which contain $a_2$ 
and $a_3$, respectively. We have $A' \in \fM_j$ 
(because $A' \cap C_j \neq \emptyset$, hence $A'$ must be contained 
in $C_j$), and $A'' \in \fM_i$ (by a similar argument). Consequently,
we have the embracings 
\begin{equation}  \label{eqn:2.15}
B' \emb A' \ \mbox{ and } \ B'' \emb A'',
\end{equation}
where $B' \in \fM_j$ and $B'' \in \fM_i$ are the blocks of $\pi$
which appear in the (Min) condition stated for $\fM_j$ and for 
$\fM_i$, respectively. But let us observe that from the hypotheses 
given on $a_1, a_2, a_3, a_4$, we can also infer that 
\begin{equation}  \label{eqn:2.16}
B'' \emb A' \ \mbox{ and } \ B' \emb A''.
\end{equation}
The embracings (\ref{eqn:2.16}) are easily verified by using the 
criterion from Proposition 2.10.1; for instance for the first of the 
two embracings we observe that
\begin{equation}  \label{eqn:2.14}
\min (B'') = \min (C_i) \leq a_1 < a_2 < a_3 \leq \max (C_i) 
=\max (B''),
\end{equation}
then we apply Proposition 2.10.1 to the situation where 
$\min (B'') < a_2 < \max ( B'' )$, with $a_2 \in A'$. (The 
equalities $\min (B'') = \min (C_i)$ and 
$\max (B'') = \max (C_i)$ appearing in (\ref{eqn:2.14}) follow
from how $C_i$ is defined in (\ref{eqn:2.13}), combined with the 
fact that $B'' \emb A$ for every $A \in \fM_i$.)

From the embracings listed in (\ref{eqn:2.15}) and (\ref{eqn:2.16}),
and by using Proposition 2.10.2, we find that we have either
$B' \emb B''$ or $B'' \emb B'$. Say for instance that the first of
these two possibilities takes place. Then we look at the embracings
$B' \emb B'' \emb A'$ with $A', B' \in \fM_j$, and we use the 
convexity condition (C) stated for $\fM_j$, to obtain that 
$B'' \in \fM_j$. Hence $B'' \in \fM_i \cap \fM_j$, which implies 
the desired conclusion that $i=j$.
(End of verification of Fact 3.)

\vspace{10pt}

It is clear that for any given integer $p$ such that 
$| \pi |_{\mathrm{out}} \leq p \leq | \pi |$, the Facts 2 and 3 
together provide us with a bijection between 
$\{ \rho \in NC(n) \mid \rho \geqgeq \pi$, $| \rho |=p \}$ and 
the collection of all CMin decompositions 
$\{ \fM_1, \ldots , \fM_p \}$ of the set of blocks of $\pi$. 
We will next observe that the latter CMin decompositions are in 
one-to-one correspondence with sets of blocks of $\pi$ which 
contain all the $| \pi |_{\mathrm{out}}$ outer blocks, plus 
$p - | \pi_{\mathrm{out}} |$ non-outer blocks. The 
precise description of this (very natural) one-to-one 
correspondence is given in the next Fact 4. Since a set of 
$p - | \pi |_{\mathrm{out}}$ non-outer blocks of $\pi$ can be 
chosen in exactly
$\left(  \begin{array}{c}  
| \pi | - | \pi |_{\mathrm{out}}  \\
   p    - | \pi |_{\mathrm{out}} 
\end{array}  \right)$ ways, the discussion of Fact 4 will 
actually conclude the proof of the proposition.

\vspace{10pt}
  
{\em Fact 4. (a) Let $\{ \fM_1, \ldots , \fM_p \}$ be a CMin 
decomposition of the set of blocks of $\pi$, and for every 
$1 \leq i \leq p$ let $B_i$ be the (uniquely determined) block
of $\pi$ which appears in the (Min) condition stated for $\fM_i$.
Then $\{ B_1, \ldots , B_p \}$ is a set of blocks of $\pi$ which 
contains all the outer blocks of $\pi$.

(b) Let $\{ B_1, \ldots , B_p \}$ be a set of blocks of $\pi$ which 
contains all the outer blocks. There exists a unique 
CMin decomposition $\{ \fM_1, \ldots , \fM_p \}$ of the set of 
blocks of $\pi$, such that $\{ B_1, \ldots , B_p \}$ is associated 
to $\{ \fM_1, \ldots , \fM_p \}$ in the way described in (a) above.}

\vspace{10pt}

The statement (a) in Fact 4 is trivial: an outer block $B$ of $\pi$ 
must satisfy $B = B_i$ for the unique $1 \leq i \leq p$ such that 
$B \in \fM_i$. The statement (b) is best understood from the 
perspective of the ``rooted forest'' framework discussed in Remark 
2.11. (It is actually immediate to translate Fact 4 into a general 
statement about rooted forests, upon suitable interpretation for 
what ``$\emb$'' and ``CMin'' should mean in that framework.) We will 
indicate how $\{ \fM_1, \ldots , \fM_p \}$ is constructed by starting 
from $\{ B_1, \ldots , B_p \}$, and we will leave it as an exercise 
to the interested reader to fill in the details of this 
graph-theoretic argument.

So suppose that we are given a block $A$ of $\pi$. The block $A$ must
be put into one of the sets of blocks $\fM_1, \ldots , \fM_p$, and 
we have to indicate the procedure for finding the index 
$i \in \{ 1, \ldots , p \}$ such that $A \in \fM_i$. We will describe
this procedure by referring to the Hasse diagram discussed in 
Remark 2.11. Let $B$ be the unique root ( = outer block) which lies 
in the same connected component of the Hasse diagram as $A$. There 
exists a unique path from $A$ to $B$ in the Hasse diagram (this 
happens because the connected components of the Hasse diagram are 
trees). Let us denote this path as $(A_0, A_1, \ldots , A_s)$, with 
$s \geq 0$, and where $A_0 = A$, $A_s = B$. Note that the path must 
intersect $\{ B_1, \ldots , B_p \}$ -- indeed, we have in any case 
that $A_s \in \{ B_1 , \ldots , B_p \}$, due to the hypothesis that 
$\{ B_1 , \ldots , B_p \}$ contains all the outer blocks. Let 
$r$ be the smallest number in $\{ 0, 1, \ldots , s \}$ such that 
$A_r \in \{ B_1 , \ldots , B_p \}$, and let $i$ be the index in 
$\{ 1, \ldots , p \}$ for which $A_r = B_i$. This $i$ is the index 
we want -- that is, the block $A$ gets placed into the set of
blocks $\fM_i$.  {\bf QED}

$\ $

{\bf 2.14 Remark.} For a given partition $\pi \in NC(n)$, the total
number of partitions $\rho \in NC(n)$ such that $\rho \geqgeq \pi$
is equal to $2^{ | \pi | - | \pi |_{\mathrm{out}} }$. This is obtained
by summing over $p$ in Equation (\ref{eqn:2.12}) of Proposition 2.13.
Or at a ``bijective'' level, one can note that in the proof of 
Proposition 2.13 the partitions $\rho \in NC(n)$ such that 
$\rho \geqgeq \pi$ end by being put into one-to-one correspondence 
with (arbitrarily chosen) sets of non-outer blocks of $\pi$.

$\ $

$\ $

\setcounter{section}{3}
\begin{center}
{\large\bf 3. \boldmath{$R, \eta, \Reta$} and 
\boldmath{$\bB$}, in the algebraic framework}
\end{center}

Throughout this section we fix a positive integer $k$ (the number of
\setcounter{equation}{0}
non-commuting indeterminates we are working with). We will deal with 
non-commutative distributions considered in an algebraic framework. 
The $R$ and $\eta$ series associated to such a distribution 
are reviewed in Definition 3.3, while $\Reta$ and $\bB$ are introduced 
in Definition 3.7.

$\ $

{\bf 3.1 Definition} {\em (non-commutative distributions).} 

$1^o$ We denote by $\ncpolk$ the algebra of non-commutative polynomials 
in $X_1, \ldots ,$
$X_k$. Thus $\ncpolk$ has a linear basis
\begin{equation}  \label{eqn:4.3011}
\{ \fatone \} \cup \{ X_{i_1} \cdots X_{i_n} \mid n \geq 1, \
1 \leq i_1, \ldots , i_n \leq k \},
\end{equation}
where the monomials in the basis are multiplied by concatenation.
When needed, $\bC \langle X_1, \ldots ,$
$X_k \rangle$ will be viewed as a $*$-algebra, with 
$*$-operation determined uniquely by the fact that each of 
$X_1, \ldots , X_k$ is selfadjoint.

$2^o$ Let $( \cM , \varphi )$ be a non-commutative probability space; that
is, $\cM$ is a unital algebra over $\bC$, and $\varphi : \cM \to \bC$ is a 
linear functional, normalized by the condition that $\varphi ( 1_{\cM} ) = 1$. 
For $x_1, \ldots , x_k \in \cM$, the {\bf joint distribution}
of $x_1, \ldots , x_k$ is the linear functional 
$\mu_{x_1, \ldots , x_k} : \ncpolk \to \bC$ which acts on the linear 
basis (\ref{eqn:4.3011}) by the formula
\begin{equation}  \label{eqn:4.302}
\left\{  \begin{array}{lcl}
\mu_{x_1, \ldots , x_k} ( \fatone ) & = & 1                        \\
\mu_{x_1, \ldots , x_k} ( X_{i_1} \cdots X_{i_n} ) & = & 
                           \varphi ( x_{i_1} \cdots x_{i_n} ),     \\
      & & \forall \, n \geq 1, \ 1 \leq i_1, \ldots , i_n \leq k. 
\end{array}   \right.
\end{equation}

$3^o$ As already mentioned in the introduction, we will denote
\begin{equation}  \label{eqn:4.303}
\Dalg (k) := \{ \mu : \ncpolk \to \bC \mid \mu \mbox{ linear, }
\mu ( \fatone ) = 1 \} .
\end{equation}
Note that, unlike in the $C^*$-context, there is no positivity requirement in the definition
of $\Dalg (k)$.
It is immediate that $\Dalg (k)$ is precisely the set of linear functionals 
on $\ncpolk$ which can appear as joint distribution for some $k$-tuple
$x_1, \ldots , x_k$ in a non-commutative probability space.

$\ $

{\bf 3.2 Definition} {\em (series and their coefficients).}

$1^o$ Recall from the introduction that 
$\bC_0 \langle \langle z_1, \ldots , z_k \rangle \rangle$ 
denotes the space of power series with complex coefficients and with 
vanishing constant term, in $k$ non-commuting indeterminates 
$z_1, \ldots , z_k$. The general form of a series 
$f \in \bC_0 \langle \langle z_1, \ldots , z_k \rangle \rangle$ is thus
\begin{equation}  \label{eqn:5.101}
f( z_1, \ldots , z_k) = \sum_{n=1}^{\infty} \
\sum_{i_1, \ldots , i_n =1}^k  \ \alpha_{(i_1, \ldots , i_n)}
z_{i_1} \cdots z_{i_n},
\end{equation}
where the coefficients  $\alpha_{(i_1, \ldots , i_n)}$ are from $\bC$. 

$2^o$ For $n \geq 1$ and $1 \leq i_1, \ldots , i_n \leq k$ we will 
denote by
\[
\cf_{(i_1, \ldots , i_n)} :
\bC_0 \langle \langle z_1, \ldots , z_k \rangle \rangle \to \bC
\]
the linear functional which extracts the coefficient of 
$z_{i_1} \cdots z_{i_n}$ in a series 
$f \in \bC_0 \langle \langle z_1, \ldots , z_k \rangle \rangle$.
Thus for $f$ written as in Equation (\ref{eqn:5.101}) we have 
$\cf_{(i_1, \ldots , i_n)} (f) = \alpha_{(i_1, \ldots , i_n)}$.

$3^o$ Suppose we are given a positive integer $n$, some indices 
$i_1, \ldots , i_n \in \{ 1, \ldots , k \}$, and a partition 
$\pi \in NC(n)$. We define a (generally non-linear) functional 
\[
\cf_{(i_1, \ldots , i_n) ; \pi} :
\bC_0 \langle \langle z_1, \ldots , z_k \rangle \rangle \to \bC ,
\]
as follows. For every block $B = \{ b_1, \ldots , b_m \}$ of $\pi$, 
with $1 \leq b_1 < \cdots < b_m \leq n$, let us use the notation 
\[
(i_1, \ldots , i_n) \vert B \ := \ (i_{b_1}, \ldots , i_{b_m})
\in \{ 1, \ldots , k \}^m.
\]
Then we define
\begin{equation}  \label{eqn:5.102}
\cf_{(i_1, \ldots , i_n); \pi} (f) \ := \
\prod_{B \ \mathrm{block \ of} \ \pi} \ \cf_{(i_1, \ldots , i_n)|B} (f),
\ \ \forall \, f \in \ncserk .
\end{equation}
(For example if we had $n=5$ and $\pi = \{ \{ 1,4,5 \} , \{ 2,3 \} \}$, 
and if $i_1, \ldots , i_5$ would be some fixed indices from 
$\{ 1, \ldots , k \}$, then the above formula would become
\[
\cf_{(i_1, i_2, i_3, i_4, i_5) ; \pi } (f) \ = \
\cf_{(i_1, i_4, i_5)} (f) \cdot
\cf_{(i_2, i_3)} (f),
\]
$f \in \bC_0 \langle \langle z_1, \ldots , z_k \rangle \rangle$.)

$\ $

{\bf 3.3 Definition} {\em (the series $M, R, \eta$).}

Let $\mu$ be a distribution in $\Dalg (k)$. We will work with three series 
$M_{\mu}, R_{\mu}, \eta_{\mu} \in$
$\bC_0 \langle \langle z_1, \ldots , z_k \rangle \rangle$ that are 
associated to $\mu$, and are defined as follows.

$1^o$ The {\bf moment series} of $\mu$ will be denoted by $M_{\mu}$. Its
coefficients are defined by 
\[
\cf_{(i_1, \ldots , i_n)} (M_{\mu}) = 
\mu (X_{i_1} \cdots X_{i_n} ), \ \ \forall \,
n \geq 1, \ \forall \, 1 \leq i_1, \ldots , i_n \leq k.
\]

$2^o$ The {\bf R-transform} of $\mu$ will be denoted by $R_{\mu}$. The 
coefficients of $R_{\mu}$ are defined by a formula which expresses them as 
polynomial expressions in the coefficients of $M_{\mu}$:
\begin{equation}  \label{eqn:4.2}
\cf_{(i_1, \ldots , i_n)} (R_{\mu}) = \sum_{\pi \in NC(n)} 
s( \pi ) \cdot \cf_{(i_1, \ldots , i_n); \pi} (M_{\mu}), 
\ \ \forall \, n \geq 1, \ \forall \,
1 \leq i_1, \ldots , i_n \leq k,
\end{equation}
where on the right-hand side of (\ref{eqn:4.2}) we used the notation 
for generalized coefficients from Definition 3.2.3, and where 
$\{ s( \pi ) \mid \pi \in \cup_{n=1}^{\infty} NC(n) \}$ is a special
family of coefficients (not depending on $\mu$). For a given 
$\pi \in NC(n)$, the precise description of $s( \pi )$ goes as follows: 
consider the Kreweras complement 
$K( \pi ) = \{ B_1, \ldots , B_p \} \in NC(n)$, and define 
$s( \pi ) := s_{|B_1|} \cdots s_{|B_p|}$, where the $s_m$ are signed 
Catalan numbers,  $s_m = (-1)^{m-1} (2m-2)!/ ( (m-1)!m!)$ for $m \geq 1$. 

The explicit description of the coefficients $s( \pi )$ is probably less 
illuminating than explaining that they appear in the following way.
The Equations in (\ref{eqn:4.2}) are equivalent to another family of 
equations of the same form, where the roles of $M_{\mu}$ and $R_{\mu}$
are switched (that is, the coefficients of $M_{\mu}$ are written as 
polynomials expressions in the coefficients of $R_{\mu}$). The $s( \pi )$
are chosen such that in this equivalent family of equations we only have
plain summations:
\begin{equation}  \label{eqn:4.203}
\cf_{(i_1, \ldots , i_n)} (M_{\mu}) = \sum_{\pi \in NC(n)} 
\cf_{(i_1, \ldots , i_n); \pi} (R_{\mu}), 
\ \ \forall \, n \geq 1, \ \forall \,
1 \leq i_1, \ldots , i_n \leq k.
\end{equation}

$3^o$ The {\bf \boldmath{$\eta$}-series} of $\mu$ will be denoted by 
$\eta_{\mu}$. The procedure for defining $\eta_{\mu}$ in terms of $M_{\mu}$ 
is analogous to the one used for defining $R_{\mu}$, only that now we are 
using the set $\Int (n)$ of interval partitions instead of $NC(n)$. The 
precise formula giving the coefficients of $\eta_{\mu}$ is
\begin{equation}  \label{eqn:4.4}
\cf_{(i_1, \ldots , i_n)} ( \eta_{\mu} ) = \sum_{\pi \in \Int (n)} 
(-1)^{1+ | \pi |} \cf_{(i_1, \ldots , i_n); \pi} ( M_{\mu} ), 
\ \ \forall \, n \geq 1, \ \forall \, 1 \leq i_1, \ldots , i_n \leq k.
\end{equation}  
The choice of the values ``$\pm 1$'' on the right-hand side of 
(\ref{eqn:4.4}) is made so that the reverse connection between the 
coefficients of $M_{\mu}$ and of $\eta_{\mu}$ is described by plain
summations:
\begin{equation}  \label{eqn:4.204}
\cf_{(i_1, \ldots , i_n)} (M_{\mu}) = \sum_{\pi \in \Int (n)} 
\cf_{(i_1, \ldots , i_n); \pi} ( \eta_{\mu} ), \ \ \forall \, n \geq 1, 
\ \forall \, 1 \leq i_1, \ldots , i_n \leq k.
\end{equation}  

$\ $

{\bf 3.4 Remark.}
$1^o$ The formulas connecting the moment series $M_{\mu}$ to the 
series $R_{\mu}$ and $\eta_{\mu}$ are well-known, and are usually 
stated as relations between certain multi-linear functionals 
(moment functionals and cumulant functionals) on non-commutative 
probability spaces. More precisely, the fomulas for $R_{\mu}$
relate to the concept of {\em free cumulant functionals} introduced in 
\cite{S94}, while the formulas for $\eta_{\mu}$ relate to the 
{\em Boolean cumulant functionals} which go all the way back to
\cite{vW73}.

$2^o$ It is clear that $\mu \mapsto M_{\mu}$ is a bijection from $\Dalg (k)$ 
onto $\ncserk$. Since the formulas which define $R_{\mu}$ and $\eta_{\mu}$
in terms of $M_{\mu}$ are reversible (in the way explained in the parts 
$2^o$ and $3^o$ of the above definition), it is immediate that
$\mu \mapsto R_{\mu}$ and $\mu \mapsto \eta_{\mu}$ also are bijections 
from $\Dalg (k)$ onto $\ncserk$; these are the bijections displayed in 
(\ref{eqn:1.10}) and (\ref{eqn:1.105}) of the introduction. 

$3^o$ Let $\mu$ be a distribution in $\Dalg (k)$. The Equation
(\ref{eqn:4.203}) describing the passage from $R_{\mu}$ to $M_{\mu}$ has
a straightforward extension to a summation formula which gives the
{\em generalized} coefficients of $M_{\mu}$ in terms of those of $R_{\mu}$.
This formula is
\begin{equation}  \label{eqn:4.30}
\cf_{(i_1, \ldots , i_n); \rho} ( M_{\mu} )
 = \sum_{ \begin{array}{c}
{\scriptstyle \pi \in NC(n),}  \\
{\scriptstyle \pi \leq \rho} 
\end{array}  } \
\cf_{(i_1, \ldots , i_n); \pi} ( R_{\mu} ), 
\end{equation}
holding for any $\rho \in NC(n)$ (and where the original formula 
(\ref{eqn:4.203}) corresponds to the case when $\rho = 1_n$).

A similar statement holds in connection to the passage from $\eta_{\mu}$
to $M_{\mu}$ -- one obtains a summation formula which gives the 
generalized coefficients of $M_{\mu}$ in terms of those of $\eta_{\mu}$,
extending Equation (\ref{eqn:4.204}). 

$4^o$ In the analytic theory of distributions of 1 variable, 
the definition of the $\eta$-series of a probability measure $\mu$ on 
$\bR$ appears usually as ``$\eta = \Psi / ( 1+ \Psi )$'',
where $\Psi$ is defined by an integral formula and corresponds, in the 
case when $\mu$ has compact support, to the moment series of $\mu$
-- see for instance the presentation at the beginning of Section 2 of
\cite{BB05}. The next proposition shows that such an approach can be also 
used in our multi-variable setting.

$\ $

{\bf 3.5 Proposition.} Let $\mu$ be a distribution in $\Dalg (k)$. We have
\begin{equation}  \label{eqn:5.01}
\eta_{\mu} = M_{\mu} /(1+M_{\mu}), 
\end{equation}
where the division on the right-hand side of (\ref{eqn:5.01}) stands for the 
commuting product $M_{\mu} (1+M_{\mu})^{-1}$ in the ring $\ncserk$. 
Conversely, $M_{\mu}$ can be obtained from $\eta_{\mu}$ by the formula
\begin{equation}  \label{eqn:5.02}
M_{\mu} = \eta_{\mu} / (1 - \eta_{\mu} ).
\end{equation}

$\ $

{\bf Proof.} We will verify the relation 
\begin{equation}  \label{eqn:4.70}
M_{\mu} = \eta_{\mu} + M_{\mu} \cdot \eta_{\mu} ,
\end{equation}
out of which (\ref{eqn:5.01}) and (\ref{eqn:5.02}) follow via easy algebraic 
manipulations.

We will fix for the whole proof some integers $n \geq 1$ and 
$1 \leq i_1, \ldots , i_n \leq k$, and we will verify the 
equality of the coefficients of $z_{i_1} z_{i_2} \cdots z_{i_n}$ in the 
series on the two sides of (\ref{eqn:4.70}). Our computations will rely on
the immediate observation that
\[
\Int (n) = \{ 1_n \} \cup \Int^{(1)}(n) \cup \cdots \cup \Int^{(n-1)}(n),
\]
disjoint union, where for $1 \leq m \leq n-1$ we denote
\[
\Int^{(m)} (n) := \{ \pi \in \Int (n) \mid \{ m+1, \ldots , n \} 
                     \mbox{ is a block of } \pi \} .
\]
We will also use the obvious fact that for every $1 \leq m \leq n-1$ we have 
a natural bijection $\Int^{(m)} (n) \ni \pi \mapsto \pi ' \in \Int (m)$, 
where $\pi '$ is obtained from $\pi$ by removing its right-most block 
$\{ m+1 , \ldots , n \}$. 

So then, compute:
\[
\cf_{(i_1, \ldots , i_n) } (M_{\mu}) \ = \
\sum_{\pi \in \Int (n)} 
\cf_{ (i_1, \ldots , i_n); \pi } ( \eta_{\mu} ) \
\mbox{ (by Equation (\ref{eqn:4.204}))}
\]
\[
= \cf_{ (i_1, \ldots , i_n); 1_n } ( \eta_{\mu} )
+ \sum_{m=1}^{n-1} \Bigl( \sum_{\pi \in \Int^{(m)} (n)} 
\cf_{ (i_1, \ldots , i_n); \pi } ( \eta_{\mu} ) \Bigr)
\]
\[
= \cf_{ (i_1, \ldots , i_n) } ( \eta_{\mu} ) + \sum_{m=1}^{n-1} 
\Bigl( \sum_{\pi ' \in \Int (m)} \cf_{ (i_1, \ldots , i_m); \pi ' } 
    ( \eta_{\mu} ) \Bigr)
\cdot \cf_{ (i_{m+1}, \ldots , i_n) } ( \eta_{\mu} ) 
\]
\[
= \cf_{ (i_1, \ldots , i_n) } ( \eta_{\mu} ) + \sum_{m=1}^{n-1} 
\cf_{ (i_1, \ldots , i_m) } (M_{\mu})  
\cdot \cf_{ (i_{m+1}, \ldots , i_n) } ( \eta_{\mu} )  \ 
\mbox{ (by Equation (\ref{eqn:4.204}))}
\]
\[
= \cf_{ (i_1, \ldots , i_n) } ( \eta_{\mu} + M_{\mu} \cdot \eta_{\mu} ),
\]
as required.  {\bf QED}

$\ $

{\bf 3.6 Remark.} Since our presentation in this section emphasizes the 
parallelism between $R$ and $\eta$, let us briefly mention that there 
exists a counterpart for Equation (\ref{eqn:5.01}) in the theory of 
the $R$-transform -- but this is a more complicated, implicit equation
involving $M_{\mu}$ and $R_{\mu}$. This latter equation is not used in 
the present paper (for a presentation of how it looks and how it is 
derived, we refer to \cite{NS06}, Lecture 16).

$\ $

{\bf 3.7 Definition.} Refer to the bijections $R$ and $\eta$ from 
$\Dalg (k)$ onto $\ncserk$ which were observed in Remark 3.4.2. We define
two new bijections:
\begin{equation}   \label{eqn:5.1}
\bB := R^{-1} \circ \eta : \Dalg (k) \to \Dalg (k),
\end{equation}
and
\begin{equation}  \label{eqn:5.2}
\Reta := \eta \circ R^{-1} : \ncserk \to \ncserk .
\end{equation}

$\ $

{\bf 3.8 Remark.}
$1^o$ We have now formally defined all the bijections which appear in
the commutative diagram (\ref{eqn:1.11}) from the introduction.
(The commutativity of this diagram
is ensured by the very definition of $\bB$ and $\Reta$.) 

$2^o$ As explained in the introduction, $\bB$ stands for 
``Boolean Bercovici--Pata bijection'', while $\Reta$ gets its name 
from the formula
\begin{equation}  \label{eqn:5.3}
\Reta (R_{\mu} ) = \eta_{\mu}, \ \ \forall \, \mu \in \Dalg (k)
\end{equation}
(it is the transformation of $\ncserk$ which ``converts $R$ to $\eta$'').

$3^o$ We will next prove that $\Reta$ can also be described by an 
explicit formula via summations over non-crossing partitions. This will 
be the same formula as indicated in Equation (\ref{eqn:1.8}) of Theorem 1,
with the difference that we will now state and prove the formula for an 
arbitrary $f \in \ncserk$, rather than just for series in the smaller set
$\cR_c (k) \subseteq \ncserk$. We mention that in the particular case 
when $k=1$, some formulas equivalent to (\ref{eqn:5.11}) and 
(\ref{eqn:5.12}) of the next proposition have appeared in \cite{L03}
(cf Equation (4.10) and the proof of Equation (5.1) in that paper).

$\ $

{\bf 3.9 Proposition.} Let $f,g$ be series in $\ncserk$ such that
$\Reta (f) =g$. Then:

$1^o$ For every $n \geq 1$ and $1 \leq i_1, \ldots , i_n \leq k$ we have
\begin{equation}  \label{eqn:5.11}
\cf_{(i_1, \ldots , i_n)} (g) 
 = \sum_{ \begin{array}{c}
{\scriptstyle  \pi \in NC(n),}  \\
{\scriptstyle \pi \leqleq 1_n}
\end{array}  } \ \cf_{(i_1, \ldots , i_n); \pi} (f).
\end{equation}

$2^o$ For every $n \geq 1$ and $1 \leq i_1, \ldots , i_n \leq k$ we have
\begin{equation}  \label{eqn:5.12}
\cf_{(i_1, \ldots , i_n)} (f) 
 = \sum_{ \begin{array}{c}
{\scriptstyle  \rho \in NC(n),}  \\
{\scriptstyle \rho \leqleq 1_n}
\end{array}  } \ (-1)^{1+| \rho |} \cf_{(i_1, \ldots , i_n); \rho} (g).
\end{equation}

$\ $

{\bf Proof.} By the definition of $\Reta$, there exists a 
distribution $\mu \in \Dalg (k)$ such that $R_{\mu} = f$, $\eta_{\mu} = g$. 

$1^o$ We calculate as follows:
\[
\cf_{(i_1, \ldots , i_n)} (g) =
\cf_{(i_1, \ldots , i_n)} ( \eta_{\mu} )
\]
\[
= \sum_{\rho \in \Int(n)} (-1)^{1+| \rho |} 
\cf_{(i_1, \ldots , i_n); \rho } ( M_{\mu} )
\ \mbox{ (by Equation (\ref{eqn:4.4})) }
\]
\[
= \ \sum_{\rho \in \Int(n)} (-1)^{1+| \rho |} \Bigl( 
\sum_{   \begin{array}{c}
{\scriptstyle \pi \in NC(n),} \\
{\scriptstyle \pi \leq \rho} 
\end{array}  } \ 
\cf_{(i_1, \ldots , i_n); \pi } ( R_{\mu} )   \Bigr)
\ \mbox{ (by Equation (\ref{eqn:4.30})) }
\]
\begin{equation} \label{eqn:4.8}
= \ \sum_{\pi \in NC(n)} \Bigl( 
\sum_{   \begin{array}{c}
{\scriptstyle \rho \in \Int(n),} \\
{\scriptstyle \rho \geq \pi} 
\end{array}  } \ (-1)^{1+| \rho |} \Bigr) 
\cf_{(i_1, \ldots , i_n); \pi } ( f ) ,
\end{equation} 
where at the last equality sign we performed a change in the order of 
summation. 

Now let us fix for the moment a partition $\pi \in NC(n)$. We claim that 
\begin{equation} \label{eqn:4.9}
\sum_{
\begin{array}{c}
{\scriptstyle \rho \in \Int(n),} \\
{\scriptstyle \rho \geq \pi} 
\end{array}  } \ (-1)^{1+| \rho |}  \ = \ \left\{
\begin{array}{ll}
1 & \mbox{ if $\pi \leqleq 1_n$}  \\
0 & \mbox{otherwise.}
\end{array}  \right.
\end{equation} 
Indeed, the set $\{ \rho \in \Int (n) \ | \ \rho \geq \pi \}$ has a 
smallest element $\widetilde{\rho}$, with $| \widetilde{\rho} |$ =
$| \rho |_{\mathrm{out}}  =: r$ (cf Remark 2.9.3). Hence we have
\[
\sum_{ \begin{array}{c}
{\scriptstyle \rho \in \Int(n),} \\
{\scriptstyle \rho \geq \pi} 
\end{array}  } \ (-1)^{1+| \rho |}  
= \sum_{ \begin{array}{c}
{\scriptstyle \rho \in \Int(n),} \\
{\scriptstyle \rho \geq \widetilde{\rho} } 
\end{array}  } \ (-1)^{1+| \rho |}  
= \sum_{m=1}^r
\left(  \begin{array}{c}  r-1 \\ m-1  \end{array} \right)
\cdot (-1)^{m+1} 
\]
\[
 = \ (1-1)^{r-1} \ = \
\left\{   \begin{array}{ll}
1 & \mbox{ if $| \rho |_{\mathrm{out}} =1$}  \\
0 & \mbox{otherwise}
\end{array}  \right\} \ = \
\left\{   \begin{array}{ll}
1 & \mbox{ if $\pi \leqleq 1_n$}  \\
0 & \mbox{otherwise,}
\end{array}  \right.
\]
as required, where during the calculation we used the immediate fact that
for every $1 \leq m \leq r$ there are 
$\left(  \begin{array}{c}  r-1 \\ m-1  \end{array} \right)$
partitions $\rho \in \Int (n)$ such that $\rho \geq \widetilde{\rho}$ and 
$| \rho | = m$.

The formula (\ref{eqn:5.11}) now follows, when (\ref{eqn:4.9}) is 
substituted in (\ref{eqn:4.8}).

$2^o$ Before starting on the calculation which leads to (\ref{eqn:5.12}),
let us note that it is straightforward to extend the Equation
(\ref{eqn:5.11}) proved in $1^o$ to a formula expressing a 
{\em generalized} coefficient of $g$ in terms of the generalized 
coefficients of $f$. (This is analogous to how (\ref{eqn:4.203}) was 
extended to (\ref{eqn:4.30}) in Remark 3.4.3.) The precise formula extending 
(\ref{eqn:5.11}) is 
\begin{equation} \label{eqn:4.60}
\cf_{(i_1, \ldots , i_n); \rho} (g) = \sum_{  \begin{array}{c}
{\scriptstyle  \pi \in NC(n),}  \\
{\scriptstyle  \pi \leqleq \rho}  
\end{array} } \ \
\cf_{(i_1, \ldots , i_n); \pi} (f),
\end{equation}
holding for an arbitrary $\rho \in NC(n)$, and where the original Equation 
(\ref{eqn:5.11}) corresponds to the case $\rho = 1_n$.

We now start from the right-hand side of (\ref{eqn:5.12}), and substitute 
$\cf_{(i_1, \ldots , i_n); \rho} (g)$ in terms of generalized coefficients 
of $f$, as indicated by Equation (\ref{eqn:4.60}). We get:
\[
\sum_{  \begin{array}{c}
{\scriptstyle  \rho \in NC(n),}  \\
{\scriptstyle  \rho \leqleq 1_n}  
\end{array}  } \  (-1)^{ 1 + | \rho | }  \Bigl(
\sum_{ \begin{array}{c}
{\scriptstyle  \pi \in NC(n),}  \\
{\scriptstyle  \pi \leqleq \rho}  
\end{array} } \  
\cf_{(i_1, \ldots , i_n); \pi} (f)  \Bigr) ;
\]
after changing the order of summation over $\pi$ and $\rho$, this becomes
\begin{equation}  \label{eqn:4.10}
\sum_{  \begin{array}{c}
{\scriptstyle  \pi \in NC(n),}  \\
{\scriptstyle  \pi \leqleq 1_n}  
\end{array} } \ \Bigl(
\sum_{  \begin{array}{c}
{\scriptstyle  \rho \in NC(n)}  \\
{\scriptstyle  such \ that }    \\
{\scriptstyle  \pi \leqleq \rho \leqleq 1_n}  
\end{array}  } \  (-1)^{1+ | \rho |} \Bigr) \cdot
\cf_{(i_1, \ldots , i_n); \pi} (f) .
\end{equation}
The sum over $\rho$ which appears in (\ref{eqn:4.10}) 
can be treated exactly as we did with (\ref{eqn:4.9}), but where this 
time instead of a counting argument in $\Int (n)$ we now invoke 
Proposition 2.13. The reader should have no difficulty to verify 
that what we get is the following: for a fixed $\pi \in NC(n)$ such 
that $\pi \leqleq 1_n$, we have
\[
\sum_{  \begin{array}{c}
{\scriptstyle  \rho \in NC(n)}  \\
{\scriptstyle  such \ that }  \\
{\scriptstyle  \pi \leqleq \rho \leqleq 1_n}  
\end{array}  } \  (-1)^{1+ | \rho |}  \ = \ 
\left\{ \begin{array}{cl}
1  &  \mbox{ if $\pi = 1_n$}  \\
0  &  \mbox{ otherwise.}
\end{array}  \right.
\]
Substituting this in (\ref{eqn:4.10}) leads to the conclusion that the 
expression considered there is equal to 
$\cf_{(i_1, \ldots , i_n)} (f)$, as required.  {\bf QED}

$\ $

$\ $

\setcounter{section}{4}
\begin{center}
{\large\bf 4. \boldmath{$\cD_c (k)$}, and infinite divisibility with 
respect to \boldmath{$\freeplus$} and \boldmath{$\Uplus$} }
\end{center}

In this section, $k$ is a fixed positive integer.
\setcounter{equation}{0}

$\ $

{\bf 4.1 Definition.} $1^o$ Refer to the space $\Dalg (k)$ introduced 
in Definition 3.1.3. We denote
\begin{equation}   \label{eqn:6.1}
\cD_c (k) = \left\{  \mu \in \Dalg (k)  \begin{array}{cc}
\vline  &  \exists \mbox{ $C^*$-probability space $( \cM , \varphi )$} \\
\vline  &  \mbox{and selfadjoint elements $x_1, \ldots , x_k \in \cM$} \\
\vline  &  \mbox{such that $\mu_{x_1, \ldots , x_k} = \mu$}
\end{array}  \right\} 
\end{equation}
(where the fact that $( \cM , \varphi )$ is a $C^*$-probability space
means that $\cM$ is a unital $C^*$-algebra, and that $\varphi : \cM \to \bC$ 
is a positive linear functional such that $\varphi ( 1_{\cM} ) = 1$). 

The notation ``$\cD_c (k)$'' is chosen to remind of
``distributions with compact support'' -- indeed, in the case when $k=1$
we have a natural identification between $\cD_c (1)$ and the set of 
probability distributions with compact support on $\bR$.

$2^o$ Refer to the bijections $R, \eta : \Dalg (k) \to \ncserk$ from 
Remark 3.4.2. We will denote 
\begin{equation}   \label{eqn:6.2}
\cR_c (k) := \{ f \in \ncserk \mid \exists \, \mu \in \cD_c (k) 
\mbox{ such that } R_{\mu} = f \} 
\end{equation}
and respectively
\begin{equation}   \label{eqn:6.3}
\cE_c (k) := \{ f \in \ncserk \mid \exists \, \mu \in \cD_c (k) 
\mbox{ such that } \eta_{\mu} = f \} .
\end{equation}
We thus have bijections
\[
\cD_c (k) \ni \mu \mapsto R_{\mu} \in \cR_c (k) \mbox{ and }
\cD_c (k) \ni \mu \mapsto \eta_{\mu} \in \cE_c (k),
\]
as indicated in (\ref{eqn:1.2}) and (\ref{eqn:1.5}) of the introduction
section.

$\ $

Clearly, the set of distributions $\cD_c (k)$ is much smaller than 
$\Dalg (k)$. In what follows we will use the characterization of 
$\cD_c (k)$ given by the next proposition. This characterization is 
most likely a ``folklore'' fact; for the reader's convenience, we include 
an outline of the proof.

$\ $

{\bf 4.2 Proposition.} Let $\mu$ be a functional in $\Dalg (k)$. Then 
$\mu$ belongs to $\cD_c (k)$ if and only if it satisfies the following 
two conditions:

(i) $\mu (P^*P) \geq 0$, $\forall \, P \in \ncpolk$.

(ii) There exists a constant $\gamma >0$ such that 
\begin{equation}  \label{eqn:6.4}
| \mu (X_{i_1} \cdots X_{i_n} ) | \leq \gamma^n, \ \ \
\forall \, n \geq 1, \ \forall \, 1 \leq i_1, \ldots , i_n \leq k.
\end{equation}

$\ $

{\bf Proof.} The necessity of the conditions (i) and (ii) is immediate, and 
left as exercise. We will outline the argument for their sufficiency. 
The proof is of course a variation of the GNS construction, the only special
point that has to be addressed is the boundedness of the left multiplication 
operators. 

So suppose that $\mu$ satisfies (i) and (ii). The positivity condition (i)
allows us to create a Hilbert space $\cH$ and a linear map
$\ncpolk \ni P \mapsto \widehat{P} \in \cH$, such that the image of 
this map is a dense subspace of $\cH$, and such that the inner product on 
$\cH$ is determined by the formula
\[
\langle \widehat{P} , \widehat{Q} \rangle = \mu (Q^*P), \ \ \forall
\, P,Q \in \ncpolk .
\]

By using the boundedness condition (ii) we will prove the inequality
\begin{equation}  \label{eqn:6.6}
|| \widehat{X_iP} || \leq \gamma || \widehat{P} ||, \ \ \forall \,
1 \leq i \leq k, \ \forall \, P \in \ncpolk ,
\end{equation}
where $\gamma > 0$ is the constant appearing in (ii).
This amounts to proving that 
\begin{equation}  \label{eqn:6.7}
\mu ( P^* X_i^2 P ) \leq \gamma^2 \mu ( P^*P ), \ \ \forall \,
1 \leq i \leq k, \ \forall \, P \in \ncpolk .
\end{equation}
We will obtain (\ref{eqn:6.7}) by a repeated application of the 
Cauchy--Schwarz inequality in $\cH$, which says that 
\begin{equation}  \label{eqn:6.8}
| \mu ( P^* Q ) | \leq \mu ( P^*P )^{1/2} \mu (Q^*Q )^{1/2}, 
\ \ \forall \, P,Q \in \ncpolk .
\end{equation}
So let us fix $i \in \{ 1, \ldots , k \}$ and $P \in \ncpolk$ and let 
us use (\ref{eqn:6.8}) with $Q= X_i^2 P$. We obtain
\begin{equation}  \label{eqn:6.9}
\mu ( P^* X_i^2 P ) \leq \mu ( P^*P )^{1/2} \mu (P^* X_i^4 P )^{1/2}. 
\end{equation}
Then let us use (\ref{eqn:6.8}) once again, but this time with 
$Q = X_i^4 P$. We get (after taking both sides to power $1/2$) that
\[
\mu ( P^* X_i^4 P )^{1/2}  \leq \mu ( P^*P )^{1/4} 
\mu (P^* X_i^8 P )^{1/4}; 
\]
and replacing the latter inequality into (\ref{eqn:6.9}) leads to
\begin{equation}  \label{eqn:6.10}
\mu ( P^* X_i^2 P ) \leq \mu ( P^*P )^{3/4} \mu (P^* X_i^8 P )^{1/4}. 
\end{equation}
It is immediate how this trick can be iterated (use Cauchy--Schwarz with 
$Q= X_i^8P$, then with $Q= X_i^{16}P$, etc), to obtain that 
\begin{equation}  \label{eqn:6.11}
\mu ( P^* X_i^2 P ) \leq \mu ( P^*P )^{(2^n-1)/2^n} 
\mu (P^* X_i^{2^{n+1}} P )^{1/2^n}, \ \ \forall \, n \geq 1. 
\end{equation}

We can now use the condition (ii) to get an upper bound on the factor 
$\mu (P^* X_i^{2^{n+1}} P )^{1/2^n}$ on the right-hand side of the
inequality (\ref{eqn:6.11}). Indeed,  let us write 
$P = \sum_{j=1}^m \alpha_j P_j$,
where for every $1 \leq j \leq m$ we have that $\alpha_j$ is in 
$\bC \setminus \{ 0 \}$ and that $P_j$ a monomial of length $l_j$ in the 
variables $X_1, \ldots , X_k$. Then
\[
P^* X_i^{2^{n+1}} P = \sum_{j,j' =1}^m 
\overline{\alpha_j} \alpha_{j'} (P_j^* X_i^{2^{n+1}} P_{j'}),
\]
which implies that
\begin{align*}
\mu ( P^* X_i^{2^{n+1}} P ) 
& \leq  \sum_{j,j' =1}^m 
| \alpha_j \alpha_{j'} | \cdot | \mu (P_j^* X_i^{2^{n+1}} P_{j'}) |  \\
& \leq  \sum_{j,j' =1}^m 
| \alpha_j \alpha_{j'} | \cdot \gamma^{l_j + 2^{n+1} + l_{j'}}
\mbox{  (by condition (ii)) }                                         \\
& = \gamma^{2^{n+1}} C,
\end{align*}
where 
$C := \sum_{j,j' =1}^m | \alpha_j \alpha_{j'} | \cdot \gamma^{l_j + l_{j'}}$
is a constant which depends only on $P$ (but not on $n$). We thus obtain 
that the right-hand side of (\ref{eqn:6.11}) is bounded from above by 
\[
\mu (P^*P)^{(2^n -1)/2^n} \cdot \Bigl( \gamma^{2^{n+1}} C \Bigr)^{1/2^n},
\]
and (\ref{eqn:6.7}) follows when we let $n \to \infty$.

Finally, by using (\ref{eqn:6.6}) it is immediatly seen that one can define 
a family of bounded linear operators $T_1, \ldots , T_k \in B( \cH )$,
determined by the formula
\[
T_i \widehat{P} \  = \ \widehat{X_iP}, \ \ \forall \, 1 \leq i \leq k, 
\ \, \forall \, P \in \ncpolk.
\]
We leave it as an easy exercise to the reader to check that 
$T_i =T_i^*$, $1 \leq i \leq k$, and that if one considers the 
$C^*$-probability space $( B( \cH ), \varphi )$ with 
$\varphi (T) := \langle T \widehat{1}, \widehat{1} \rangle$,
$T \in B( \cH )$, then the joint distribution of $T_1, \ldots , T_k$
in $( B ( \cH ), \varphi )$ is precisely the functional $\mu$ that 
we started with. It follows that $\mu \in \cD_c (k)$, as required.
{\bf QED}

$\ $

We now come to the operations $\freeplus$ and $\Uplus$ that were mentioned 
in the introduction. It will be convenient to consider them in the larger 
algebraic framework provided by the space $\Dalg (k)$. Each of these 
two operations has its own theory, developped in connection to a form of 
independence for non-commutative random variables. We will briefly
comment on this in the Remark 4.4 below, but we will not need to go into 
details about non-commutative independence. Indeed, for the approach used
in this paper (where the $R$-transform and $\eta$-series play the main
role), we can simply regard $\freeplus$ and $\Uplus$ as the binary operations 
that are linearized by $R$ and $\eta$, respectively.

$\ $

{\bf 4.3 Definition.} Let $\mu$ and $\nu$ be distributions in $\Dalg (k)$.

$1^o$ The {\bf free additive convolution} $\mu \freeplus \nu$ is the 
unique distribution in $\Dalg (k)$ which has $R$-transform equal to
\begin{equation}   \label{eqn:40.11}
R_{\mu \freeplus \nu} = R_{\mu} + R_{\nu}.
\end{equation}

$2^o$ The {\bf Boolean convolution} $\mu \Uplus \nu$ is the unique 
distribution in $\Dalg (k)$ which has $\eta$-series equal to
\begin{equation}   \label{eqn:40.12}
\eta_{\mu \Uplus \nu} = \eta_{\mu} + \eta_{\nu}.
\end{equation}

$\ $

{\bf 4.4 Remark.}
The above definition reverses the order of how things are usually 
considered in the literature -- usually $\freeplus$ and $\Uplus$ are 
considered first, and then $R$ and $\eta$ appear as linearizing transforms
for these two operations. The way how $\freeplus$ and $\Uplus$ are usually 
considered is in connection to the concepts of {\em free independence} and 
respectively {\em Boolean independence} for subsets of a non-commutative 
probability space $( \cM , \varphi )$. More precisely, suppose that we have
elements $x_1, \ldots , x_k, y_1, \ldots ,$ $y_k \in \cM$ such that the 
joint distribution of $x_1, \ldots , x_k$ is equal to $\mu$, and the joint 
distribution of $y_1, \ldots , y_k$ is equal to $\nu$. If 
$\{ x_1, \ldots , x_k \}$ is freely independent from 
$\{ y_1, \ldots , y_k \}$ in $( \cM , \varphi )$, then the joint 
distribution of the $k$-tuple $x_1 + y_1, \ldots , x_k + y_k$ is equal to
$\mu \freeplus \nu$; while if $\{ x_1, \ldots , x_k \}$ is Boolean 
independent from $\{ y_1, \ldots , y_k \}$ in $( \cM , \varphi )$, then 
the joint distribution of the $k$-tuple $x_1 + y_1, \ldots , x_k + y_k$ is 
equal to $\mu \Uplus \nu$. 

In this paper we will not need to review the precise definitions of free
and of Boolean independence. We need however to mention one fact about 
$\freeplus$ and $\Uplus$ which comes out of the approach via non-commutative
independence, namely that:
\begin{equation}   \label{eqn:6.12}
\mu , \nu \in \cD_c (k) \ \Longrightarrow \ 
\mu \freeplus \nu , \
\mu \Uplus \nu  \in \cD_c (k).
\end{equation}
Indeed, if $\mu , \nu \in \cD_c (k)$, then it can be shown that 
$x_1, \ldots , x_k,y_1, \ldots , y_k$ from the preceding paragraph can 
always be found to be selfadjoint elements in a $C^*$-probability space. 
Since the elements $x_1+y_1, \ldots , x_k+y_k$ will then also be selfadjoint, 
it follows that the convolutions $\mu \freeplus \nu$ and $\mu \Uplus \nu$ 
are still in $\cD_c (k)$. Thus we can (and will) also view $\freeplus$ and 
$\Uplus$ as binary operations on $\cD_c (k)$. 

Let us also record here that, as a consequence of (\ref{eqn:6.12}) and of 
the Equations (\ref{eqn:40.11}) and (\ref{eqn:40.12}) in Definition 4.3, 
it is immediate that the sets of series 
$\cR_c (k),  \cE_c (k) \subseteq \ncserk$ are closed under addition.

$\ $ 

We now move to discuss infinite divisibility. We discuss first the case of 
$\freeplus$.

$\ $

{\bf 4.5 Definition.}
$1^o$ Let $\mu$ be in $\cD_c (k)$. If for every positive integer $N$
there exists a distribution $\mu_N \in \cD_c (k)$ such that
\[
\underbrace{\mu_N \freeplus \cdots \freeplus \mu_N}_{N \ \mathrm{times}}
= \mu,
\]
then we will say that $\mu$ is
{\bf \boldmath{$\freeplus$}-infinitely divisible.}
The set of all distributions $\mu \in \cD_c (k)$ which are 
$\freeplus$-infinitely divisible will be denoted by $\Dinfdiv (k)$.

$2^o$ We denote
\begin{equation}   \label{eqn:6.15}
\Rinfdiv (k) = \{ f \in \cR_c (k) \mid f = R_{\mu} \mbox{ for a 
distribution } \mu \in \Dinfdiv (k) \}.
\end{equation}

$\ $

{\bf 4.6 Remark.} Infinite divisibility with respect to $\freeplus$ 
relates to how $\cR_c (k)$ behaves under the operation of multiplication 
by scalars from $(0, \infty )$. Let us record here that we have:
\begin{equation}   \label{eqn:6.13}
\Bigl( f \in \cR_c (k), t \in [ 1, \infty ) \Bigr) \ \Rightarrow \
tf \in \cR_c (k).
\end{equation}
This is a non-trivial fact, which appears in connection to how 
$R$-transforms behave under compressions by free projections -- see
Lecture 14 of \cite{NS06} for more details.

On the other hand, $\cR_c (k)$ is not closed under multiplication by 
scalars from $(0,1)$. For a fixed series $f \in \cR_c (k)$ we have in 
fact that
\begin{align*}
tf \in \cR_c (k), \ \forall \, t \in (0, 1) 
& \Leftrightarrow \ \frac{1}{N} f \in \cR_c (k), \ \forall \,
     N \in \bN \setminus \{ 0 \}                                  \\
& \Leftrightarrow \ f \in \Rinfdiv (k), 
\end{align*}
where the first of these equivalences follows from (\ref{eqn:6.13}), 
and the second one is a direct consequence of Definition 4.5.

$\ $

{\bf 4.7 Remark.} Following the above considerations about infinite 
divisibility for $\freeplus$, it would be now natural to do the parallel
discussion and introduce the corresponding notations for $\Uplus$. But 
it turns out that no new notations are needed, as $\cE_c (k)$ is closed 
under multiplication by scalars from $(0, \infty )$ (and consequently, all 
the distributions in $\cD_c (k)$ are $\Uplus$-infinitely divisible). This 
fact is proved in the next 
proposition, by using an operator model for how to achieve the 
multiplication of an $\eta$-series by a scalar $t \in (0,1)$. We mention 
here that in the case $k=1$ another proof of this proposition can be 
given by using complex analysis methods (specific to the case $k=1$ only);
see Theorem 3.6 of \cite{SW97}. To our knowledge, the case $k \geq 2$ was
not treated before (it is e.g mentioned as an open problem in the recent
thesis \cite{W04}).

$\ $

{\bf 4.8 Proposition.} If $f \in \cE_c (k)$ and $t \in ( 0, \infty )$, 
then $tf \in \cE_c (k)$.

$\ $

{\bf Proof.} Since we know that $\cE_c (k)$ is closed under addition,
it suffices to do the case when $t \in (0,1)$. We fix for the whole proof
a series $f \in \cE_c (k)$ and a number $t \in (0,1)$. We denote by $\mu$
the unique distribution in $\cD_c (k)$ such that $\eta_{\mu} = f$; the goal
of the proof is to find a distribution $\nu \in \cD_c (k)$ such that 
$\eta_{\nu} = tf$.

Let $x_1, \ldots , x_k$ be selfadjoint elements in a $C^*$-probability 
space $( \cM , \varphi )$ such that the joint distribution of 
$x_1, \ldots , x_k$ is equal to $\mu$. 
By considering the GNS representation of $\varphi$ we may assume, 
without loss of generality, that $\cM = B( \cH )$ for a Hilbert space
$\cH$, and that $\varphi$ is the vector-state given by a unit vector
$\xi_o \in \cH$ (that is, $\varphi (x) = \langle x \xi_o , \xi_o \rangle$
for every $x \in B( \cH )$).

Let us consider a new Hilbert space
\begin{equation}  \label{eqn:6.17}
\cK := \cV \oplus ( \cH \ominus \bC \xi_o ) \oplus ( \cH \ominus \bC \xi_o ),
\end{equation}
where $\cV$ is a Hilbert space of dimension 2, spanned by two vectors 
$\Omega_1, \Omega_2$ such that 
\begin{equation}  \label{eqn:6.18}
|| \Omega_1 || = 1 = || \Omega_2 ||, \ \ 
\langle \Omega_1 , \Omega_2 \rangle = t-1.
\end{equation}
We consider moreover two isometric operators 
$J_1, J_2 : \cH \to \cK$, defined by
\begin{equation}  \label{eqn:6.19}
\left\{  \begin{array}{cccl}
J_1 ( \alpha \xi_o + \xi ) & = & \alpha \Omega_1 \oplus \xi \oplus 0 &  \\
J_2 ( \alpha \xi_o + \xi ) & = & \alpha \Omega_2 \oplus 0 \oplus \xi , &
\forall \, \alpha \in \bC, \ \forall \, \xi \in
\cH \ominus \bC \xi_o .  
\end{array}  \right.
\end{equation}
It is immediately verified that the adjoints 
$J_1^*, J_2^* \in B( \cK , \cH )$ are described by the formula
\begin{equation}  \label{eqn:6.20}
\left\{  \begin{array}{cccl}
J_1^* ( v \oplus \xi_1 \oplus \xi_2 ) & = 
               & \langle v , \Omega_1 \rangle  \xi_o + \xi_1  &  \\
J_2^* ( v \oplus \xi_1 \oplus \xi_2 ) & = 
               & \langle v , \Omega_2 \rangle  \xi_o + \xi_2 , &
\forall \, \alpha \in \bC, \ \forall \, 
\xi_1, \xi_2 \in \cH \ominus \bC \xi_o .  
\end{array}  \right.
\end{equation}
As a consequence of (\ref{eqn:6.18})--(\ref{eqn:6.20}), we have
\begin{equation}  \label{eqn:6.21}
J_1^* J_2 = J_2^* J_1 = (t-1) P_o,
\end{equation}
where $P_o \in B( \cH )$ is the orthogonal projection onto the 1-dimensional 
space $\bC \xi_o$ (that is, $P_o \xi = \langle \xi , \xi_o \rangle \xi_o$ 
for every $\xi \in B( \cH )$).

Consider the operators $y_1, \ldots , y_k \in B( \cK )$ defined by
\begin{equation}  \label{eqn:6.22}
y_i := J_1 x_i J_1^* + J_2 x_i J_2^*, \ \ 1 \leq i \leq k.
\end{equation}
Let $\nu \in \cD_c (k)$ be the joint distribution of $y_1, \ldots , y_k$
in the $C^*$-probability space $( B( \cK ), \psi )$, where $\psi$ is the 
vector-state given on $B( \cK )$ by the unit-vector $\Omega$ defined as
\begin{equation}  \label{eqn:6.23}
\Omega := \frac{1}{\sqrt{2t}} \Bigl( \, ( \Omega_1 + \Omega_2 )
\oplus 0 \oplus 0 \, \Bigr) \in \cK .
\end{equation}
We want to obtain an explicit formula for the coefficients of the moment 
series $M_{\nu}$. So let us fix a positive integer $n$ and 
some indices $i_1, \ldots , i_n \in \{ 1, \ldots , k \}$, and let us compute:
\[
\cf_{(i_1, \ldots , i_n)} (M_{\nu}) = 
\psi( y_{i_1} \cdots y_{i_n} ) =
\langle y_{i_1} \cdots y_{i_n} \Omega , \Omega \rangle \ \ 
\mbox{ (by the def. of $\nu$ and of $\psi$)}
\]
\[
= \langle (J_1 x_{i_1} J_1^* + J_2 x_{i_1} J_2^*) \cdots 
(J_1 x_{i_n} J_1^* + J_2 x_{i_n} J_2^*) \Omega , \Omega \rangle \ \ 
\mbox{ (by the def. of $y_1, \ldots , y_k$)}        
\]
\[
= \sum_{r_1, \ldots , r_n =1}^2 \langle (J_{r_1} x_{i_1} J_{r_1}^*) \cdots 
(J_{r_n} x_{i_n} J_{r_n}^* ) \Omega , \Omega \rangle \ \ 
\]
\[
= \sum_{r_1, \ldots , r_n =1}^2 \langle x_{i_1} (J_{r_1}^*J_{r_2}) x_{i_2}
\cdots (J_{r_{n-1}}^* J_{r_n}) x_{i_n} (J_{r_n}^* \Omega ) , 
(J_{r_1}^* \Omega ) \rangle
\]
\begin{equation}  \label{eqn:6.24}
= \frac{t}{2} \sum_{r_1, \ldots , r_n =1}^2 \langle x_{i_1} 
(J_{r_1}^*J_{r_2}) x_{i_2} \cdots (J_{r_{n-1}}^* J_{r_n}) x_{i_n} \xi_o , 
\xi_o \rangle ,
\end{equation}
where at the last equality sign we took into account that (as is immediately
verified) $J_1^* \Omega = J_2^* \Omega = \sqrt{t/2} \xi_o$.
The next thing to be taken into account is that, for $r,r' \in \{ 1,2 \}$,
we have:
\[
J_r^*  J_{r'} = \left\{  \begin{array}{ccc}
1_{B( \cH )}  & \mbox{if} & r=r'   \\
(t-1)P_o      & \mbox{if} & r \neq r'.
\end{array}  \right.
\]
So for any $n$-tuple $(r_1, \ldots , r_n) \in \{ 1, 2 \}^n$, the operator
\begin{equation}  \label{eqn:6.25}
x_{i_1} (J_{r_1}^*J_{r_2}) x_{i_2} (J_{r_2}^* J_{r_3}) \cdots 
(J_{r_{n-1}}^* J_{r_n}) x_{i_n} \in B( \cH )
\end{equation}
really depends only on the set of positions $m \in \{ 1, \ldots , n-1 \}$
where we have $r_m \neq r_{m+1}$. If we write this set of positions as
$\{ m_1, \ldots , m_p \}$ with $1 \leq m_1 < m_2 < \cdots < m_p \leq n-1$,
then the product considered in (\ref{eqn:6.25}) equals
\begin{equation}  \label{eqn:6.26}
(x_{i_1} \cdots x_{i_{m_1}}) \bigl( (t-1)P_o \bigr)
(x_{i_{m_1+1}} \cdots x_{i_{m_2}}) \bigl( (t-1)P_o \bigr) \cdots
\bigl( (t-1)P_o \bigr) 
(x_{i_{m_p+1}} \cdots x_{i_n}).
\end{equation}
It is moreover immediately seen that when one applies the vector-state 
$\varphi = \langle \centerdot \ \xi_o , \xi_o \rangle$ to the operator in 
(\ref{eqn:6.26}), the result is 
\begin{equation}  \label{eqn:6.27}
(t-1)^p \langle x_{i_1} \cdots x_{i_{m_1}}  \xi_o , \xi_o \rangle 
\langle x_{i_{m_1+1}} \cdots x_{i_{m_2}} \xi_o , \xi_o \rangle  \cdots
\langle x_{i_{m_p+1}} \cdots x_{i_n} \xi_o , \xi_o \rangle .
\end{equation}
(Note: It is not ruled out that the set 
$\{ m  \mid 1 \leq m \leq n-1$, $r_m \neq r_{m+1} \}$ could be empty. 
The formula (\ref{eqn:6.27}) still holds in this case, the quantity 
appearing there being then equal to 
$\langle x_{i_1} \cdots x_{i_n} \xi_o , \xi_o \rangle$.)

It is convenient that in the calculations shown in the preceding paragraph 
we encode the sequence $1 \leq m_1 < \cdots < m_p \leq n-1$ by the interval
partition $\pi = \{ B_1, \ldots, B_{p+1} \}$ where 
$B_1 = [1, m_1 ] \cap \bZ$,
$B_2 = (m_1, m_2 ] \cap \bZ, \ldots ,B_p = (m_{p-1}, m_p ] \cap \bZ$,
$B_{p+1} = (m_p, n ] \cap \bZ$. (In the case when 
$\{ m_1, \ldots , m_p \} = \emptyset$, we take $\pi$ to be $1_n$, the 
partition with only one block.) The quantity in (\ref{eqn:6.27}) then
becomes 
\begin{equation}  \label{eqn:6.28}
(t-1)^{| \pi |-1} \cf_{(i_1, \ldots , i_n); \pi} (M_{\mu}),
\end{equation}
where the generalized coefficient
$\cf_{(i_1, \ldots , i_n); \pi} (M_{\mu})$
is considered in the sense of Definition 3.2.3.
Note moreover that for every given 
partition $\pi \in \Int (n)$ there are exactly two $n$-tuples 
$(r_1, \ldots , r_n ) \in \{ 1,2 \}^n$ for which the set 
$\{ m \mid 1 \leq m \leq n-1$, $r_m \neq r_{m+1} \}$ is encoded by $\pi$
(one of these two $n$-tuples has $r_1 =1$, and the other has $r_1 =2$).

If we now return to the expression in (\ref{eqn:6.24}), and if in that
summation formula we replace $n$-tuples 
$(r_1, \ldots , r_n ) \in \{ 1,2 \}^n$ by partitions $\pi \in \Int (n)$,
then we obtain:
\[
\cf_{(i_1, \ldots , i_n)} ( M_{\nu} ) = 
t \sum_{\pi \in \Int (n)} (t-1)^{| \pi | -1}
\cf_{(i_1, \ldots , i_n); \pi} ( M_{\mu} ).
\]

It is more suggestive to re-write the above equation in the form
\begin{equation}  \label{eqn:6.29}
\cf_{(i_1, \ldots , i_n)} \Bigl( \frac{t-1}{t} M_{\nu} 
\Bigr) = \sum_{\pi \in \Int (n)} 
\cf_{(i_1, \ldots , i_n); \pi} \Bigl( (t-1) M_{\mu} \Bigr) ,
\end{equation}
holding for every $n \geq 1$ and every $1 \leq i_1, \ldots , i_n \leq k$.
From (\ref{eqn:6.29}) we see that the series $((t-1)/t) M_{\nu}$ is 
obtained from $(t-1) M_{\mu}$ by exactly the formula expressing a moment
series in terms of the corresponding $\eta$-series -- cf Equation 
(\ref{eqn:4.204}) in Definition 3.3. But we saw in Proposition 3.5 how the 
latter formula can be written in a compressed way, in terms of the 
series themselves rather than in terms of coefficients. Applied to the
situation at hand, Proposition 3.5 will thus give us that 
\begin{equation}  \label{eqn:6.30}
\frac{t-1}{t} M_{\nu} = \frac{ (t-1)M_{\mu} }{ 1- (t-1) M_{\mu} } .
\end{equation}

Finally, we use Equation (\ref{eqn:6.30}) in order to compute $\eta_{\nu}$.
We leave it as a straightforward exercise to the reader to check that when 
we write $\eta_{\nu}$ =
$M_{\nu} / (1+ M_{\nu})$ and then replace $M_{\nu}$ in terms of $M_{\mu}$ 
by using (\ref{eqn:6.30}), what comes out is simply that
$\eta_{\nu} = ( t M_{\mu} )/ (1+ M_{\mu} )$ = 
$t \eta_{\mu}$ = $tf$, as we wanted.  {\bf QED}

$\ $

$\ $

\setcounter{section}{5}
\begin{center}
{\large\bf 5. Proofs of Theorems 1 and 1'}
\end{center}

Throughout this section, $k$ is a fixed positive integer.

$\ $

{\bf 5.1 Remark.} Theorems 1 and 1' take place in the framework of 
\setcounter{equation}{0}
$\cD_c (k)$, but in their proofs it will be nevertheless useful to rely 
on occasion on the larger algebraic framework provided by $\Dalg (k)$.
For example: when we need to construct a distribution in $\cD_c (k)$
which satisfies certain requirements, it may come in handy to first
observe a distribution $\mu \in \Dalg (k)$ which satisfies the given
requirements, and then to verify (by using Proposition 4.2) that $\mu$
belongs in fact to the subset $\cD_c (k) \subseteq \Dalg (k)$.

In connection to the above, it will be convenient to place the next 
definition (for convergence of sequences)
in the larger framework of $\Dalg (k)$.

$\ $

{\bf 5.2 Definition.} $1^o$ Let $\mu$ and $( \mu_N )_{N=1}^{\infty}$ be 
distributions in $\Dalg (k)$. The notation ``$\toN \mu_N =  \mu$'' will be 
used to mean that $( \mu_N )_{N=1}^{\infty}$ {\bf converges in moments}
to $\mu$, i.e. that
\begin{equation}  \label{eqn:7.1}
\lim_{N \to \infty} \mu_N (P) = \mu (P), \ \ \forall \, P \in \ncpolk .
\end{equation}

$2^o$ Let $f$ and $( f_N )_{N=1}^{\infty}$ be series in $\ncserk$. The 
notation ``$\toN f_N = f$'' will be used to mean that $( f_N )_{N=1}^{\infty}$ 
{\bf converges coefficientwise} to $f$, i.e. that
\begin{equation}  \label{eqn:7.2}
\lim_{N \to \infty} \cf_{(i_1, \ldots , i_n)} (f_N) =
\cf_{(i_1, \ldots , i_n)} (f), 
\ \ \forall \, n \geq 1, \ \forall \, 1 \leq i_1, \ldots , i_n \leq k.
\end{equation}

$\ $

We now start on a sequence of lemmas which will gradually build towards 
the statements of Theorems 1 and 1'.

$\ $

{\bf 5.3 Lemma.}  Let $\mu$ and $( \mu_N )_{N=1}^{\infty}$ be distributions
in $\Dalg (k)$. The following three statements are equivalent:

\begin{center}
(1) $\toN \mu_N = \mu$; \hspace{1cm}
(2) $\toN R_{\mu_N} = R_{\mu}$; \hspace{1cm}
(3) $\toN \eta_{\mu_N} = \eta_{\mu}$.
\end{center}

$\ $

{\bf Proof.} It is immediate that the convergence in moments from (1)
is equivalent to a statement (1') referring to the coefficientwise 
convergence of the corresponding moment series,

\begin{center}
(1') $\toN M_{\mu_N} =  M_{\mu}$.
\end{center}

\noindent
On the other hand,
it is immediate that we have (1') $\Leftrightarrow$ (2) and
(1') $\Leftrightarrow$ (3), due to the explicit formulas relating the 
coefficients of the series $M$, $R$, $\eta$ via (finite!) summations 
over partitions, as presented in Definition 3.3 above.  {\bf QED}

$\ $

{\bf 5.4 Lemma.}  Let $f$ and $( f_N )_{N=1}^{\infty}$ be series in
$\ncserk$, such that $\toN f_N = f$, and let $( t_N )_{N=1}^{\infty}$ be 
a sequence in $(0, \infty )$ such that 
$\lim_{N \to \infty} t_N =  \infty$. For every 
$N \geq 1$ let us consider the series 
\[
g_N := t_N \cdot \Reta \bigl( \frac{1}{t_N} f_N \bigr) \ \mbox{ and } \
h_N := t_N \cdot \Reta^{-1} \bigl( \frac{1}{t_N} f_N \bigr) 
\]
in $\ncserk$. Then $\toN g_N = f$ and $\toN h_N = f$.

$\ $

{\bf Proof.} For every $n \geq 1$ and $1 \leq i_1, \ldots , i_n \leq k$ 
we have:
\[
\cf_{(i_1, \ldots , i_n)} (g_N) = 
t_N \cdot \cf_{(i_1, \ldots , i_n)} \Bigl( \, \Reta \bigl ( \frac{1}{t_N}
f_N \bigr) \, \Bigr)
\]
\[
= \ t_N \sum_{ \begin{array}{c}
{\scriptstyle \pi \in NC(n),}  \\
{\scriptstyle \pi \leqleq 1_n} 
\end{array}  }  \ \cf_{(i_1, \ldots , i_n); \pi} \bigl ( \frac{1}{t_N}
f_N \bigr) \ \mbox{ (by Proposition 3.9.1)}
\]
\begin{equation}  \label{eqn:7.3}
= \ \sum_{ \begin{array}{c}
{\scriptstyle \pi \in NC(n),}  \\
{\scriptstyle \pi \leqleq 1_n} 
\end{array}  }  \ t_N^{1-| \pi|} \cf_{(i_1, \ldots , i_n); \pi} (f_N),
\end{equation}
where at the last equality sign we used the obvious homogeneity property
of $\cf_{(i_1, \ldots , i_n); \pi}$. When we make $N \to \infty$ in
(\ref{eqn:7.3}), the only term which survives is the one corresponding to
$\pi = 1_n$, and it follows that 
\[
\lim_{N \to \infty} \cf_{(i_1, \ldots , i_n)} (g_N) \ = \ 
\lim_{N \to \infty} \cf_{(i_1, \ldots , i_n)} (f_N) \ = \ 
\cf_{(i_1, \ldots , i_n)} (f).
\]
This proves that $\toN g_N = f$. The argument for $\toN h_N = f$ is similar, 
with the only difference that we now use Proposition 3.9.2 instead of 
3.9.1.  {\bf QED}

$\ $

{\bf 5.5 Lemma.} Let $( \mu_N )_{N=1}^{\infty}$ be in 
$\Dalg (k)$ and let $p_1 < p_2 < \cdots < p_N < \cdots$ be a sequence 
of positive integers. 

$1^o$ Suppose there exists $\mu \in \Dalg (k)$ such that 
$\toN \ \underbrace{\mu_N \Uplus \cdots \Uplus \mu_N}_{p_N} =  \mu$.
Then it follows that 
$\toN \ \underbrace{\mu_N \freeplus \cdots \freeplus \mu_N}_{p_N} =
\bB ( \mu )$ (where $\bB : \Dalg (k) \to \Dalg (k)$ is the bijection 
from Definition 3.7).

$2^o$ Suppose there exists $\nu \in \Dalg (k)$ such that 
$\toN \underbrace{\mu_N \freeplus \cdots \freeplus \mu_N}_{p_N} = \nu$.
Then it follows that 
$\toN \ \underbrace{\mu_N \Uplus \cdots \Uplus \mu_N}_{p_N}$
= $\bB^{-1} ( \nu )$.

$\ $

{\bf Proof.} $1^o$ Let us denote
\[
f_N \ :=  \ \eta_{\underbrace{\mu_N \Uplus \cdots \Uplus \mu_N}_{p_N}}
\ = \ p_N \cdot \eta_{\mu_N}, \ N \geq 1,
\]
and let $f := \eta_{\mu}$. From Lemma 5.3 and the given hypothesis
it follows that $\toN f_N = f$. Thus if we let
\[
h_N := p_N \cdot \Reta^{-1} \bigl( \frac{1}{p_N} f_N \bigr) , \ N \geq 1,
\]
then Lemma 5.4 gives us that $\toN h_N = f$ as well. But let us observe 
that, for every $N \geq 1$:
\begin{align*}
h_N &=  p_N \cdot \Reta^{-1} \bigl( \eta_{\mu_N} ) \ \mbox{ (since
$f_N = p_N \cdot \eta_{\mu_N}$)}                                          \\
&= p_N \cdot R_{\mu_N} \ \mbox{ (by definition of $\Reta$, Def. 3.7)}     \\
&= R_{\underbrace{\mu_N \freeplus \cdots \freeplus \mu_N}_{p_N}}.
\end{align*}
On the other hand we can write $f= \eta_{\mu} = R_{\bB ( \mu )}$
(by definition of $\bB$). So the convergence 
$\toN h_N = f$ amounts in fact to
\[
\toN \ R_{\underbrace{\mu_N \freeplus \cdots \freeplus \mu_N}_{p_N}} \ = \
R_{\bB ( \mu )},
\]
and the conclusion that 
$\toN \ \underbrace{\mu_N \freeplus \cdots \freeplus \mu_N}_{p_N} =
\bB ( \mu )$ follows from Lemma 5.3.

$2^o$ The proof of this statement is identical to the proof of $1^o$,
where now we switch the roles of $\freeplus$ and $\Uplus$, the roles 
of $R$ and $\eta$, and we use the other part of Lemma 5.4.  {\bf QED}

$\ $

{\bf 5.6 Lemma.} Let $\mu$ be a distribution in $\Dalg (k)$, and consider
the series $R_{\mu}, \eta_{\mu} \in \ncserk$. The following statements
are equivalent:

(1) There exists a constant $\gamma > 0$ such that 
\[
| \mu (X_{i_1} \cdots X_{i_n}) | \leq \gamma^n, \ \ \forall \, n \geq 1,
\ \forall \, 1 \leq i_1, \ldots , i_n \leq k.
\]

(2) There exists a constant $\gamma > 0$ such that 
\[
| \cf_{(i_1, \ldots , i_n)} (R_{\mu}) | \leq \gamma^n, \ \ \forall 
\, n \geq 1, \ \forall \, 1 \leq i_1, \ldots , i_n \leq k.
\]

(3) There exists a constant $\gamma > 0$ such that 
\[
| \cf_{(i_1, \ldots , i_n)} (\eta_{\mu}) | \leq \gamma^n, \ \ \forall 
\, n \geq 1, \ \forall \, 1 \leq i_1, \ldots , i_n \leq k.
\]

$\ $

{\bf Proof.} Both the equivalences (1) $\Leftrightarrow$ (2) and
(1) $\Leftrightarrow$ (3) follow from the explicit relations via 
summations over partitions which connect the coefficients of the 
series $M_{\mu}, R_{\mu}, \eta_{\mu}$, where one uses suitable bounds for 
how many terms there are in the summations, and for the size of the 
coefficients (if there are any coefficients involved). For example, 
when proving that (1) $\Rightarrow$ (2), one uses the bound
\[
\card \Bigl( NC(n) \Bigr) \ = \ \frac{(2n)!}{n!(n+1)!} \ \leq \ 4^n,
\ \ \forall \, n \geq 1,
\]
and one also uses the fact (easily proved by induction) that the constants
``$s( \pi )$'' appearing in Equation (\ref{eqn:4.2}) from Definition 3.3.2
satisfy
\[
| s( \pi ) | \ \leq \ 4^n, \ \ \forall \, n \geq 1, \ \forall \, 
\pi \in NC(n).
\]
Suppose $\gamma > 0$ is such that (1) holds. Then 
the coefficients of the moment series $M_{\mu}$ satisfy
$| \cf_{(i_1, \ldots , i_n)} (M_{\mu}) | \ \leq \ \gamma^n$, for every
$n \geq 1$ and $1 \leq i_1, \ldots , i_n \leq k$, and more generally
\[
| \cf_{(i_1, \ldots , i_n); \pi} (M_{\mu}) | \ \leq \ \gamma^n, \ \ \forall 
\, n \geq 1, \ \forall \, 1 \leq i_1, \ldots , i_n \leq k, \ \forall 
\, \pi \in NC(n).
\]
So then for every $n \geq 1$ and $1 \leq i_1, \ldots , i_n \leq k$ we have
\[
| \cf_{(i_1, \ldots , i_n)} (R_{\mu}) | = 
| \sum_{\pi \in NC(n)}
s( \pi ) \cdot \cf_{(i_1, \ldots , i_n); \pi} (M_{\mu}) | 
\leq \ \sum_{\pi \in NC(n)} 4^n \cdot \gamma^n \leq (16 \gamma )^n,
\]
and (2) holds, with $\gamma$ replaced by $16 \gamma$. The arguments for 
(2) $\Rightarrow$ (1) and for (1) $\Leftrightarrow$ (3) are similar, 
and require in fact smaller corrections for $\gamma$.  {\bf QED}

$\ $

{\bf 5.7 Lemma.} We have $\Rinfdiv (k) = \cE_c (k)$, where the subsets
$\Rinfdiv (k)$ and $\cE_c (k)$ of $\ncserk$ are as in Definition 4.5.2 and 
in Definition 4.1.2, respectively.

$\ $

{\bf Proof.} ``$\subseteq$'' Let $f$ be a series in $\Rinfdiv (k)$, about 
which we want to show that $f  \in \cE_c (k)$. Let $\mu$ be the unique 
distribution in $\Dalg (k)$ such that $\eta_{\mu} = f$; proving that 
$f \in \cE_c (k)$ is equivalent to proving that $\mu \in \cD_c (k)$.
We will prove the latter fact, by verifying that $\mu$ satisfies the 
conditions (i) and (ii) from Proposition 4.2. 

The verification of (ii) is immediate, in view of Lemma 5.6. Indeed, 
since $f \in \Rinfdiv (k)$
$\subseteq \cR_c (k)$, we know there exists a
distribution $\nu \in \cD_c (k)$ such that $R_{\nu} = f$. By using the 
condition (ii) for $\nu$ and the equivalence (1) $\Leftrightarrow$ (2)
in Lemma 5.6, we find that there exists $\gamma > 0$ such that 
$| \cf_{(i_1, \ldots , i_n)} (f) | \leq \gamma^n$ for every $n \geq 1$
and every $1 \leq i_1, \ldots , i_n \leq k$. But then we can use the 
fact that $f = \eta_{\mu}$ and the equivalence (1) $\Leftrightarrow$ (3)
in Lemma 5.6, to obtain that $\mu$ satisfies (ii).

In order to verify that $\mu$ satisfies condition (i), we proceed 
as follows. For every $N \geq 1$ let us consider the series 
\begin{equation}  \label{eqn:7.4}
g_N := N \cdot \Reta \bigl( \frac{1}{N} f \bigr) \in \ncserk,
\end{equation}
and the unique distribution $\mu_N \in \Dalg (k)$ such that 
$\eta_{\mu_N} = g_N$. We have $\toN g_N = f$, by Lemma 5.4. This convergence
can also be written as $\toN \eta_{\mu_N} = \eta_{\mu}$, and it gives 
us that $\toN \mu_N = \mu$, by Lemma 5.3. But now let us observe that for
every $N \geq 1$ we have
\[
\frac{1}{N} f \in \cR_c (k) \ \mbox{ (because $f \in \Rinfdiv (k)$)}
\]
\[
\Rightarrow \ \Reta \bigl( \frac{1}{N} f \bigr) \in \cE_c (k) \ 
\mbox{ (because $\Reta$ maps $\cR_c (k)$ onto $\cE_c (k)$)}
\]
\[
\Rightarrow \ g_N = N \cdot \Reta \bigl( \frac{1}{N} f \bigr) \in 
\cE_c (k) \ \mbox{ (because $\cE_c (k)$ is closed under addition)}
\]
\[
\Rightarrow \ \mu_N \in \cD_c (k) \ 
\mbox{ (by the definition of $\cE_c (k)$).}
\]
In particular it follows that every $\mu_N$ satisfies condition (i), and 
then it is clear that the limit in moments 
$\mu = \lim_{N \to \infty} \mu_N$ has to satisfy (i) too. 

``$\supseteq$'' Let us observe that it suffices to prove the weaker 
inclusion 
\begin{equation}  \label{eqn:7.50}
\cE_c (k)  \subseteq \cR_c (k).
\end{equation}
Indeed, if (\ref{eqn:7.50}) is known, then for an arbitrary series 
$f \in \cE_c (k)$ we get that:
\[
tf \in \cE_c (k), \ \ \forall \, t \in (0,1)
\mbox{ (by Proposition 4.8)}
\]
\[
\Rightarrow \ tf \in \cR_c (k) \ \ \forall \, t \in (0,1)
\mbox{ (by (\ref{eqn:7.50}))}
\]
\[
\Rightarrow \ f \in \Rinfdiv (k) \ \ \mbox{ (by Remark 4.6).}
\]

Hence for this part of the proof it suffices if we fix a series 
$f \in \cE_c (k)$, and prove that $f \in \cR_c (k)$. The argument for
this is pretty much identical to the one shown above, in the proof of
the inclusion $\subseteq$. That is, we consider the unique distribution
$\nu \in \Dalg (k)$ such that $R_{\nu} = f$, and we prove that 
$\nu \in \cD_c (k)$, by verifying that it satisfies the conditions (i)
and (ii) from Proposition 4.2. The verification of (ii) proceeds exactly
as in the proof of $\subseteq$ (we look at the distribution 
$\mu \in \cD_c (k)$ which has $\eta_{\mu} = f = R_{\nu}$, and we use 
Lemma 5.6 twice, in connection to $\mu$, $f$ and $\nu$). The verification
of (i) also proceeds on the same lines as shown in the proof of 
$\subseteq$, with the difference that instead of the series $g_N$ from
(\ref{eqn:7.4}) we now look at 
\begin{equation}  \label{eqn:7.5}
h_N := N \cdot \Reta^{-1} \bigl( \frac{1}{N} f \bigr) , \ \ N \geq 1,
\end{equation}
and we consider the distributions $( \nu_N )_{N=1}^{\infty}$ in 
$\Dalg (k)$ which have $R_{\nu_N} = h_N$, $N \geq 1$. We leave it as an 
exercise to the reader to adjust the argument shown in the proof of
$\subseteq$ in order to verify that $\nu_N \in \cD_c (k)$ for every
$N \geq 1$, and that $\toN \nu_N = \nu$. The property (i) for $\nu$ is 
therefore obtained by passing to the limit the property (i) for the 
$\nu_N$.  {\bf QED}

$\ $

{\bf 5.8 Remark} {\em (proofs of Theorems 1 and 1').}

At this moment we are in fact only left to observe that all the statements
made in Theorems 1 and 1' are covered by the arguments shown above, as 
follows.

(a) Part $1^o$ of Theorem 1 is covered by Lemma 5.7.

(b) For part $2^o$ of Theorem 1, we observe that by its very definition
(and by the definitions of $\cR_c (k)$ and $\cE_c (k)$), the bijection
$\Reta : \ncserk \to \ncserk$ from Definition 3.7 maps 
$\cR_c (k)$ onto $\cE_c (k)$. Since $\cE_c (k) = \Rinfdiv (k)$, we get
that $\Reta$ is indeed a bijection between $\cR_c (k)$ and $\Rinfdiv (k)$.
The explicit formula given for $\Reta$ in Equation (\ref{eqn:1.8}) was 
proved in Proposition 3.9.

(c) The map $\bB : \cD_c (k) \to \Dinfdiv (k)$ in part $3^o$ of Theorem 1 is
the composition of the bijections $\eta : \cD_c (k) \to \Rinfdiv (k)$ and 
$R^{-1} : \Rinfdiv (k) \to \Dinfdiv (k)$, and is therefore itself bijective.
It is also clear that $\bB : \cD_c (k) \to \Dinfdiv (k)$ is the restriction
of the bijection $\bB : \Dalg (k) \to \Dalg (k)$ from Definition 3.7.

(d) Theorem 1' follows immediately from Lemma 5.5 and the fact, recorded 
above, that the algebraic bijection $\bB : \Dalg (k) \to \Dalg (k)$ maps 
$\cD_c (k)$ onto $\Dinfdiv (k)$.

(e) The last thing left is the compatibility (stated in part $3^o$ of 
Theorem 1) of $\bB$ with the corresponding bijection from \cite{BP99}.
This is an immediate consequence of Theorem 1', which was in fact given as 
a multi-variable counterpart for the corresponding statement in
\cite{BP99}.

$\ $

{\bf 5.9 Remark.} We conclude this section by mentioning a few facts that 
did not appear in Theorems 1 and 1' or in their proofs, but may be of
relevance for other developments related to these theorems.

$1^o$ Let $t$ be a number in $(0,1)$, and let $\mu, \mu' , \mu ''$ be 
distributions in $\Dalg (k)$ such that:

\begin{center}
(i) $R_{\mu '} = \frac{1}{1-t} R_{\mu}$, 
\hspace{.5cm} and \hspace{.5cm} 
(ii) $\eta_{\mu ''} = t \eta_{\mu '}$.
\end{center}

\noindent
Then a direct computation using the relations between $R$-transforms and 
$\eta$-series yields
\begin{equation}   \label{eqn:7.7}
R_{\mu ''} = \frac{t}{1-t} \eta_{\mu}.
\end{equation}
Let us observe moreover that if $\mu \in \cD_c (k)$ then 
$\mu '$ and $\mu ''$ belong to $\cD_c (k)$ as well; this is because 
$\cR_c (k)$ is closed under multiplication by $1/(1-t)$, and $\cE_c (k)$ 
is closed under multiplication by $t$. In the case when $t= 1/2$, these 
observations can be used to give an alternative proof for the inclusion 
(\ref{eqn:7.50}) in the proof of Lemma 5.7. 

$2^o$ A positivity phenomenon which was observed in preceding work on
multi-variable $\freeplus$-infinite divisibility is the following. Let 
$f$ be a series in $\Rinfdiv (k)$ and let $\mu \in \Dalg (k)$ be the 
distribution determined by the formula
\[
\mu ( X_{i_1} \cdots X_{i_n} ) = \cf_{(i_1, \ldots , i_n)} (f), 
\ \ \forall \, n \geq 1, \ \forall \, 1 \leq i_1 , \ldots , i_n \leq k.
\]
Then $\mu (P^*P) \geq 0$ for every polynomial $P \in \ncpolk$ which has
no constant term. Moreover, this positivity can be used to construct 
realizations of $\freeplus$-infinitely divisible distributions by operators on 
the full Fock space. For more details on this, see Sections 4.5 and 4.7 of 
\cite{S98}.

$3^o$ A property of distributions $\mu \in \cD_c (k)$ which is often 
considered is {\em traciality} ($\mu (PQ) = \mu (QP)$ for every 
$P,Q \in \ncpolk$). This did not appear in Theorems 1 and 1', and in fact 
the multi-variable Boolean Bercovici-Pata bijection does not preserve
traciality. From a combinatorial perspective, the cause of this fact is 
that the natural action of the cyclic group $\bZ_n$ on partitions of 
$\{ 1, \ldots , n \}$ does not leave invariant the set $\Int (n)$ of
interval partitions.

$\ $

\setcounter{section}{6}
\begin{center}
{\large\bf 6. A special property of the partial order 
\boldmath{$\leqleq$} }
\end{center}

The goal of this section is to prove a combinatorial result which 
\setcounter{equation}{0}
lies at the heart of the proofs of Theorems 2 and 2', and which will 
be stated precisely in Corollary 6.11. The proof of this result 
will use the various facts about non-crossing partitions that were 
reviewed in Section 2, tailored to the special situation of 
parity-preserving partitions $\theta \in NC(2n)$.

$\ $

{\bf 6.1 Remark.} 
Let $n$ be a positive integer. We refer to the above Section 2.2 for 
the definition of what it means that a partition $\theta \in NC(2n)$ 
is parity-preserving, and for the fact that 
such a partition can always be uniquely presented in the form 
$\theta = \piodd \UU \roeven$ where $\pi, \rho \in NC(n)$ are such 
that $\rho \leq K( \pi )$.

A useful remark is that for $\theta = \piodd \UU \roeven$ as above
we always have
\begin{equation}
| \theta | = |\pi| + |\rho| \geq |\pi| + |K( \pi )| = n+1,
\end{equation}
with the equality $|\theta| = n+1$ holding if and only if 
$\rho = K( \pi )$.

If $\theta \in NC(2n)$ is parity-preserving, then the blocks $A$ of 
$\theta$ such that $A \subseteq \{ 1, 3, \ldots , 2n-1 \}$ will be 
called {\bf odd blocks}, while the blocks $B$ of $\theta$ such that 
$B \subseteq \{ 2, 4, \ldots , 2n \}$ will be called {\bf even blocks}. 

Let us observe that a parity-preserving partition $\theta \in NC(2n)$
always has at least two outer blocks. Indeed, the odd block $M$ such 
that $M \ni 1$ and the even block $N$ such that $N \ni 2n$ are distinct,
and both have to be outer. If these $M$ and $N$ are the only outer 
blocks of $\theta$, then we will say (naturally) that $\theta$ 
{\bf has exactly two outer blocks}. In view of Remark 2.9.2, a 
necessary and sufficient condition for this to happen is that
\begin{equation}   \label{eqn:8.2}
\min (N) = \max (M)+1.
\end{equation}
If $\theta$ is written as $\piodd \UU \roeven$ with 
$\pi, \rho \in NC(n)$ such that $\rho \leq K( \pi )$, then it is 
immediately seen that the condition (\ref{eqn:8.2}) amounts to 
\begin{equation}   \label{eqn:8.3}
\min (N_0) = \max (M_0),
\end{equation}
where $N_0$ is the block of $\rho$ such that $N_0 \ni n$ and 
$M_0$ is the block of $\pi$ such that $M_0 \ni 1$. The condition
(\ref{eqn:8.3}) is nicely expressed in terms of the permutations 
$P_{\pi}, P_{\rho}$ associated to $\pi$ and $\rho$ in Remark 2.1. 
Indeed, it is immediate that $\min (N_0) = P_{\rho} (n)$ and
$\max (M_0) = P_{\pi}^{-1} (1)$,
so in the end we arrive to the following equivalence: for 
$\pi , \rho \in NC(n)$ such that $\rho \leq K( \pi )$ we have that
\begin{equation}   \label{eqn:8.4}
\left( \begin{array}{c}
\piodd \UU \roeven \mbox{ has exactly}  \\
\mbox{two outer blocks}
\end{array}  \right) \ \Leftrightarrow \
P_{\rho} (n) = P_{\pi}^{-1} (1).
\end{equation}

$\ $

{\bf 6.2 Remark.} Let $\theta$ be a parity-preserving partition in 
$NC(2n)$, and let us consider the Hasse diagram for the embracing 
partial order on the set of blocks of $\theta$. This is exactly as 
in Remark 2.11 of Section 2, with the additional ingredient that
the vertices of the Hasse diagram are now bicoloured (indeed,
a vertex of the Hasse diagram is a block of $\theta$, and is of one 
of the colours ``even'' or ``odd''). We warn the reader that, in 
general, this bicolouring does not turn the Hasse diagram into a 
so-called ``bipartite graph'' (i.e. it is not precluded that an 
edge of the Hasse diagram connects two vertices of the same colour); 
this issue is discussed in more detail in Remark 6.5 and in Lemmas
6.6 and 6.8 below.

In what follows, instead of talking explicitly about edges of the 
Hasse diagram, we will prefer to use the related concept of 
``parenthood'' for the vertices of the diagram. This concept is defined
in the general framework of a rooted forest, and goes as follows.
Let $A$ be a vertex of a rooted forest, and suppose that $A$ is not a root 
of the forest. Let $B$ be the unique root which lies in the same 
connected component of the forest as $A$, and consider the unique
path $(A_0, A_1, \ldots , A_s)$ from $A$ to $B$ in the forest
(with $s \geq 1$, $A_0 =A$, $A_s =B$). The vertex $A_1$ of this path
is then called the {\em parent} of the vertex $A$.

In the case of the specific rooted forest that we are dealing with here
(Hasse diagram for embracing partial order on blocks), the concept of 
parenthood can of course be also given directly in terms of embracings
of blocks. We leave it as an immediate exercise to the reader to 
check that, when proceeding on this line, the definition comes out 
as follows.

$\ $

{\bf 6.3 Definition.} Let $\theta$ be a parity-preserving partition 
in $NC(2n)$, and let $A$ be a block of $\theta$, such that $A$ is not 
outer. Then there exists a block $P$ of $\theta$, uniquely determined,
with the following two properties:

\begin{center}
(i) $P \emb A$, $P \neq A$, \hspace{.5cm} and \hspace{.5cm}
(ii) If $A'\emb A$ and $A' \neq A$ then $A' \emb P$.
\end{center}

\noindent
This block $P$ will be called the {\bf parent of $A$ with respect 
to $\theta$}, and will be denoted as $\Parent_{\theta} (A)$.

$\ $

The next proposition states a few basic properties of the 
parenthood relation defined above. The verifications of 
these properties are immediate, and are left as exercise to the 
reader.

$\ $

{\bf 6.4 Proposition.} Let $\theta$ be a parity-preserving partition 
in $NC(2n)$.

$1^o$ Let $A$ be a block of $\theta$ such that $\max (A) < 2n$. Let 
$B$ be the block of $\theta$ such that $\max (A) +1 \in B$, and 
suppose that $\max (A)+1$ is not the minimal element of $B$. Then 
$A$ is not outer, and $\Parent_{\theta} (A) =B$.

$2^o$ Let $A,B$ be blocks of $\theta$ such that 
$\min(B) = \max (A) +1$. Then either both $A$ and $B$ are outer 
blocks, or none of them is, and in the latter case we have that 
$\Parent_{\theta} (A) = \Parent_{\theta} (B)$.

$3^o$ Let $A$ be a block of $\theta$ which is not outer, denote
$\Parent_{\theta} (A) =: B$, and suppose that $A$ and $B$ have the 
same parity. Let $\theta '$ be the partition of 
$\{ 1, \ldots , 2n \}$ which is obtained from $\theta$ by joining 
together the blocks $A$ and $B$. Then $\theta '$ is a parity-preserving
partition in $NC(2n)$, and we have $\theta \leqleq \theta '$.

$\ $

{\bf 6.5 Remark.} Let $\theta$ be a parity-preserving partition in 
$NC(2n)$, and consider the Hasse diagram of the embracing partial 
order on blocks of $\theta$. We will next look at the question of 
when is it possible that an edge of the Hasse diagram connects two 
blocks of the same parity. It is
clear that if two blocks of $\theta$ are connected by an edge of the 
Hasse diagram, then there is one of the two blocks which is the 
parent of the other; due to this fact, 
the question stated above amounts to asking when is it 
possible that the blocks $A$ and $\mbox{Parent}_{\theta} (A)$ have the 
same parity (where $A$ is a non-outer block of $\theta$). This will be 
addressed in the Lemmas 6.6 and 6.8 below. We will actually be only 
interested in the situation when $\theta$ has exactly two outer blocks.

$\ $

{\bf 6.6 Lemma.} Let $\pi$ be a partition in $NC(n)$, and consider 
the parity-preserving partition 
$\theta := \piodd \UU \kpeven \in NC(2n)$. 

$1^o$ $\theta$ has exactly two outer blocks.

$2^o$ If $A$ is a block of $\theta$ which is not outer, then the block
$\Parent_{\theta} (A)$ has parity opposite from the parity of $A$.

$\ $

{\bf Proof.} $1^o$ In view of the above equivalence (\ref{eqn:8.4}), 
it suffices to observe that 
$P_{K( \pi )} (n) = P_{\pi}^{-1} (P_{1_n} (n)) = P_{\pi}^{-1} (1)$
(where we used the formula for $P_{K(\pi)}$ given in 
Equation (\ref{eqn:2.1}) of Remark 2.3).

$2^o$ We will present the argument in the case when $A$ is an odd 
block of $\theta$. (The case when $A$ is an even block is treated 
similarly, and is left as exercise.) We have 
$A= \{ 2a-1 \mid a \in A_o \}$ where $A_o$ is
a block of $\pi$. Let us denote $\min (A_o) = a'$ and $\max (A_o) = a''$. 
We have $a' > 1$ (from $a'=1$ it would follow that $A \ni 1$, and $A$
would be an outer block of $\theta$). Observe that
\[
P_{K( \pi )} (a'-1) = P_{\pi}^{-1} (P_{1_n} (a'-1)) 
                    = P_{\pi}^{-1} (a') = a''.
\]
This shows that $a'-1$ and $a''$ belong to the same block $B_o$ of 
$K( \pi )$. The block $B = \{ 2b \mid b \in B_o \}$ of $\theta$ will 
then contain the elements $2(a'-1) = \min (A) -1$ and $2a'' = \max (A)+1$. 
From Proposition 6.4.1 we infer that $\Parent_{\theta} (A) = B$, and 
in particular it follows that $\Parent_{\theta} (A)$ has parity 
opposite from the one of $A$, as required.  {\bf QED}

$\ $

{\bf 6.7 Remark.}
Let us look again at the discussion from Remark 6.2, but where now 
$\theta$ is taken to be of the special kind from the preceding lemma,
$\theta := \piodd \UU \kpeven$ for a partition $\pi \in NC(n)$. So 
we look at the Hasse diagram for the embracing partial order on 
blocks of $\theta$, and we note that in this case the Hasse diagram 
{\em is} bipartite -- indeed, Lemma 6.6 assures us that every edge 
of the Hasse diagram connects an even block of $\theta$ with an odd
block of $\theta$.

The discussion around the Hasse diagram for this special $\theta$ 
can be pushed a bit further, in order to establish a connection with 
the theory of graphs on surfaces. We make here a brief comment on 
how this goes. (The comment will not be used in what follows, but 
may be illuminating for a reader who is familiar with the theory of
graphs on surfaces, and prefers to translate various facts about 
$NC(n)$ into the language of that theory.) So let us note another 
feature of our Hasse diagram: it has exactly two 
connected components, rooted at the two outer blocks of $\theta$,
and where one of the roots (the outer block containing 1) is odd
while the other (the outer block containing $2n$) is even. Let us 
add to the Hasse diagram an edge which connects the two roots. The 
graph we obtain is a tree with $n+1$ vertices (and $n$ edges), where 
the vertices are bicoloured, and where every edge connects two 
vertices of different colours; moreover, this tree has a ``marked''
edge, namely the edge that was added to the Hasse diagram. This 
kind of structure (a bicoloured bipartite tree with one marked edge)
plays an important role in the theory of graphs on surfaces, where 
among other things it is used to produce a special factorization
``of genus zero'' of the cycle 
$1 \mapsto 2 \mapsto \cdots \mapsto n \mapsto 1$; see for instance
section 1.5 of the monograph \cite{LZ04}. In particular, the 
bicoloured tree we just encountered (by starting from 
$\theta := \piodd \UU \kpeven$, and by adding an edge to the 
corresponding Hasse diagram) can be used to retrieve the 
factorization ``$P_{\pi} P_{K( \pi )} = P_{1_n}$'' reviewed in the 
Remark 2.3 -- hence this bicoloured tree completely determines the 
partition $\pi$ we started with.

The next lemma is in some sense a converse of Lemma 6.6.

$\ $

{\bf 6.8 Lemma.} Let $\pi$ and $\rho$ be partitions in $NC(n)$ such 
that $\rho \leq K( \pi)$, and consider the parity-preserving partition 
$\theta := \piodd \UU \roeven \in NC(2n)$. Consider the following 
properties that $\theta$ may have:
\begin{tabular}[t]{cl}
(i)  &  $\theta$ has exactly two outer blocks.                \\
(ii) & If $A$ is a block of $\theta$ which is not outer, 
                                               then the block \\
     & $\Parent_{\theta} (A)$ has parity opposite from the 
                                              parity of $A$.
\end{tabular}

\noindent
If both (i) and (ii) hold, then $\rho = K( \pi )$.

$\ $

{\bf Proof.} Let $M$ and $N$ be the blocks of $\theta$ such that 
$M \ni 1$ and $N \ni 2n$. Hypothesis (i) implies that
$\min (N) = 1+ \max (M)$ (cf. Equation (\ref{eqn:8.2}) in 
Remark 6.1). We denote
\[
m := \min (N)/2 = (1+ \max (M))/2 \in \{ 1, \ldots , n \} .
\]

Let $\fN$ be the set of blocks of $\theta$ which are not outer, and 
define a function $F: \fN \to \{ 1, \ldots , n \}$ by the formula:
\[
F(X) := \left\{  \begin{array}{ll}
(1+ \max (X))/2  &  \mbox{if $X$ is an odd block}   \\
\min (X)/2       &  \mbox{if $X$ is an even block.}  
\end{array}   \right.
\]
It is immediate that $F(X) \neq m$, $\forall \, X \in \fN$, 
hence that $\card ( \, F( \fN ) \, ) \leq n-1$.

We claim that $F(X) \neq F(Y)$ whenever $X,Y \in \fN$ are 
blocks of opposite parities. Indeed, assume for instance that $X$ is odd,
$Y$ is even, and $F(X) = F(Y) =: i$. Then $\min (Y) = 2i = 1+ \max (X)$,
and Proposition 6.4.2 implies that 
$\Parent_{\theta}(X) = \Parent_{\theta} (Y)$. But this is not possible,
since (by hypothesis (ii)) the block
$\Parent_{\theta}(X)$ is even, while $\Parent_{\theta} (Y)$ is odd.

It then follows that the function $F$ is one-to-one. Indeed, if 
$X,Y \in \fN$ are such that $F(X) = F(Y)$ then $X$ and $Y$ must have 
the same parity (by the claim proved in the preceding paragraph), and 
it immediately follows that $X \cap Y \neq \emptyset$, hence that $X=Y$.

From the injectivity of $F$ we conclude that 
$\card ( \fN ) = \card ( \, F( \fN ) \, ) \leq n-1$,
and that the total number of blocks of $\theta$ is 
$|\theta | = 2+ \card ( \fN ) \leq n+1$. But it was noticed 
in Remark 6.1 that 
$\theta = \piodd \UU \roeven$ has $| \theta | \geq n+1$, with
equality if and only if $\rho = K( \pi )$. The conclusion of the lemma 
immediately follows.  {\bf QED}

$\ $

{\bf 6.9 Lemma.} Let $\pi$ be in $NC(n)$, and let $\theta$ be a 
parity-preserving partition in $NC(2n)$ such that 
$\theta \leqleq \piodd \UU \kpeven$. Let $X,Y$ be blocks of $\theta$
such that $Y = \Parent_{\theta} (X)$, and suppose that $Y$ has the same
parity as $X$. Let $\theta '$ be the partition of $\{ 1, \ldots , 2n \}$
which is obtained from $\theta$ by joining together the blocks $X$ and $Y$.
Then $\theta '$ also is a parity-preserving partition in $NC(2n)$ such 
that $\theta ' \leqleq \piodd \UU \kpeven$. 

$\ $

{\bf Proof.} We will write the argument in the case when $X$ and $Y$ are
even blocks of $\theta$. The case when $X$ and $Y$ are odd blocks is 
similar, and is left as exercise.

The fact that $\theta '$ is a parity-preserving partition in $NC(2n)$ 
follows from Proposition 6.4.3, the point of the proof is to verify 
that $\theta ' \leqleq \piodd \UU \kpeven$.

Let us write $\theta = \sigmaodd \UU \taueven$ with 
$\sigma , \tau \in NC(n)$. Since we assummed that $X,Y$ are even blocks 
of $\theta$, we thus have $X = \{ 2a \mid a \in A \}$ and
$Y = \{ 2b \mid b \in B \}$ with $A,B$ blocks of $\tau$. It is clear 
that $\theta ' = \sigmaodd \UU ( \tau ' )^{\mathrm{(even)}}$, where 
$\tau ' \in NC(n)$ is obtained from $\tau$ by joining together its 
blocks $A$ and $B$.

Consider the partition $\phi := \piodd \UU \taueven$ of 
$\{ 1, \ldots , 2n \}$; that is, $\phi$ is obtained by putting together 
the even blocks of $\theta$ with the odd blocks of
$\piodd \UU \kpeven$. Note that $\phi \in NC(2n)$; this is because we have
\begin{equation}  \label{eqn:6.50}
\sigmaodd \UU \taueven = \theta \leqleq \piodd \UU \kpeven ,
\end{equation}
which implies in particular that $\tau \leq K ( \pi )$. It is clear that
$X$ and $Y$ are blocks of $\phi$. We observe the following facts.

\vspace{10pt}

{\em Fact 1. $Y$ is the block-parent of $X$ with respect to $\phi$.}

{\em Verification of Fact 1.} Assume by contradiction that there exists
a block $Z$ of $\phi$ such that $Z \neq X,Y$ and $Y \emb Z \emb X$. The 
block $Z$ cannot be even -- in that case it would be a block of $\theta$,
and the embracings $Y \emb Z \emb X$ would contradict the hypothesis that
${\mbox{Parent}}_{\theta} (X) = Y$. So $Z$ is odd, and consequently it is
a block of $\piodd \UU \kpeven$. Since $\theta \leqleq \piodd \UU \kpeven$,
there exists a block $W$ of $\theta$ which has $\min (W) = \min (Z)$ and 
$\max (W) = \max (Z)$. But then the embracings $Y \emb Z \emb X$ are 
equivalent to $Y \emb W \emb X$; since $W \neq X,Y$ (which happens because
$W$ is odd, while $X,Y$ are even), we again obtain a 
contradiction with the hypothesis that ${\mbox{Parent}}_{\theta} (X) = Y$. 
(End of verification of Fact 1.)

\vspace{10pt}

{\em Fact 2. $\tau ' \leq K ( \pi )$. }

{\em Verification of Fact 2.} Let $\phi '$ be the partition obtained from 
$\phi$ by joining together its blocks $X$ and $Y$. It is clear that $\phi '$ 
= $\piodd \UU ( \tau ' )^{\mathrm{(even)}}$; thus proving the 
inequality $\tau ' \leq K( \pi )$ is equivalent to proving that 
$\phi ' \in NC(2n)$ (see the equivalence (\ref{eqn:2.01}) in Remark 2.2). 
But $\phi '$ is indeed in $NC(2n)$, as we see by invoking the above Fact 1
and Proposition 6.4.3. (End of verification of Fact 2.)

\vspace{10pt}

We conclude the argument as follows: look at the inequalities
$\sigma \leqleq \pi$ and $\tau \leqleq K( \pi )$ which are implied by
(\ref{eqn:6.50}), and combine the second of these inequalities with 
Fact 2, in order to get that $\tau ' \leqleq K( \pi )$. (Indeed, from 
$\tau \leq \tau ' \leq K( \pi )$ and $\tau \leqleq K( \pi )$ we get that
$\tau ' \leqleq K( \pi )$ -- this follows directly from how $\leqleq$ 
is defined.) So we have $\sigma \leqleq \pi$ and 
$\tau ' \leqleq K( \pi )$, which entails that
$\theta ' = \sigmaodd \UU ( \tau ' )^{\mathrm{(even)}}  \leqleq
\piodd \UU \kpeven$, as required.  {\bf QED}

$\ $

{\bf 6.10 Proposition.} Let $\theta$ be a parity-preserving partition 
in $NC(2n)$, which has exactly two outer blocks. There exists a unique 
partition $\pi \in NC(n)$ such that $\theta \leqleq \piodd \UU \kpeven$. 

$\ $

{\bf Proof.} Let us denote
\[
\cT := \{ \theta ' \in NC(2n) \mid \theta' \mbox{ parity-preserving and }
          \theta ' \geqgeq \theta \} .
\]
Observe that every $\theta ' \in \cT$ has exactly two outer blocks 
(because $\theta$ is like that, and by the equivalence 
(\ref{eqn:2.110}) in Remark 2.12).

Let $\widetilde{\theta} \in \cT$ be an element which is maximal with
respect to the partial order $\leqleq$; that is, $\widetilde{\theta}$
is such that if $\theta ' \in \cT$ and 
$\theta ' \geqgeq \widetilde{\theta}$, then 
$\theta ' = \widetilde{\theta}$. 

Let $X$ be a block of $\widetilde{\theta}$ which is not outer, and let
us denote $Y := \Parent_{\widetilde{\theta}} (X)$. We claim that $X$ and
$Y$ have opposite parities. Indeed, in the opposite case the partition
$\theta '$ obtained from $\widetilde{\theta}$ by joining the blocks 
$X$ and $Y$ together would still be in $\cT$, and would satisfy 
$\theta ' \geqgeq \widetilde{\theta}$, 
$\theta ' \neq  \widetilde{\theta}$, contradicting the maximality of 
$\widetilde{\theta}$. 

We thus see that $\widetilde{\theta}$ satisfies the hypotheses of 
Lemma 6.8, and must therefore be of the form $\piodd \UU \kpeven$ for
some $\pi \in NC(n)$. This proves the existence part of the lemma.

For the uniqueness part, suppose that $\pi ,\rho \in NC(n)$ are such 
that $\theta \leqleq \piodd \UU \kpeven$ and
$\theta \leqleq \rhodd \UU \kroeven$. We then consider the set
\[
\cS := \left\{ \theta ' \in NC(2n) \begin{array}{cl}
\vline &  \theta' \mbox{ parity-preserving, }
          \theta ' \geqgeq \theta ,                               \\
\vline &  \theta ' \leqleq \piodd \UU \kpeven \mbox{ and }
          \theta ' \leqleq \rhodd \UU \kroeven
\end{array} \right\} ,
\]
and we let $\widehat{\theta}$ be a maximal element in $( \cS , \leqleq )$. 

Let $X$ be a block of $\widehat{\theta}$ which is not outer, and let
us denote $Y := \Parent_{\widehat{\theta}} (X)$. We claim that $X$ and
$Y$ have opposite parities. Indeed, in the opposite case the partition
$\theta '$ obtained from $\widehat{\theta}$ by joining the blocks 
$X$ and $Y$ together would still be in $\cS$ -- where the inequalities 
$\theta ' \leqleq \piodd \UU \kpeven$ and 
$\theta ' \leqleq \rhodd \UU \kroeven$ follow from Lemma 6.9. This 
partition $\theta '$ would satisfy $\theta ' \geqgeq \widehat{\theta}$ 
and $\theta ' \neq  \widehat{\theta}$, contradicting the maximality of 
$\widehat{\theta}$. 

We thus see that $\widehat{\theta}$ satisfies the hypotheses of 
Lemma 6.8, and must therefore be of the form 
$\sigma^{\mathrm{(odd)}} \UU K( \sigma )^{\mathrm{(even)}}$ for some 
$\sigma \in NC(n)$. But then the inequality
$\widehat{\theta} \leqleq \piodd \UU \kpeven$ implies that
$\sigma \leq \pi$ and $K( \sigma ) \leq K( \pi )$, which in turn imply
that $\sigma = \pi$ (as $K(\sigma) \leq K(\pi)$ is equivalent to
$\sigma \geq \pi$). A similar argument shows that $\sigma = \rho$, and 
the desired equality $\pi = \rho$ follows.  {\bf QED}

$\ $

{\bf 6.11 Corollary.} Let $\sigma , \tau$ be two partitions
in $NC(n)$. The following two statements are equivalent:

(1) There exists $\pi \in NC(n)$ such that $\sigma \leqleq \pi$ and 
$\tau \leqleq K( \pi )$.

(2) $\tau \leq K( \sigma )$ and the associated permutations 
$P_{\sigma}$ and $P_{\tau}$ are such that
$P_{\sigma}^{-1} (1) = P_{\tau} (n)$. 

\noindent
Moreover, in the case when the statements (1) and (2) are true, the 
partition $\pi$ with the properties stated in (1) is uniquely determined.

$\ $

{\bf Proof.} ``(1)$\Rightarrow$(2)'' Fix $\pi \in NC(n)$ such that
$\sigma \leqleq \pi$ and $\tau \leqleq K( \pi )$. We have in particular 
that $\sigma \leq \pi$ and $\tau \leq K( \pi )$. So 
$\tau \leq K( \pi ) \leq K( \sigma)$, and the inequality 
$\tau \leq K( \sigma )$ follows. On the other hand it is immediate that 
we have the implications
\begin{center}
$\sigma \leqleq \pi \Rightarrow P_{\sigma}^{-1} (1) = P_{\pi}^{-1} (1)$
\hspace{.5cm} and \hspace{.5cm}
$\tau \leqleq K( \pi ) \Rightarrow P_{\tau} (n) = P_{K( \pi )} (n)$.
\end{center}
So the relation $P_{K( \pi )} (n) = P_{\pi}^{-1} (1)$ observed in the 
proof of Lemma 6.6.1 implies that $P_{\tau} (n) = P_{\sigma}^{-1} (1)$.

``(2)$\Rightarrow$(1)'' Consider the partition 
$\theta = \sigmaodd \UU \taueven$ of $\{ 1, \ldots , 2n \}$. The 
hypotheses in (2) give us that $\theta \in NC(2n)$ and that it has 
exactly two outer blocks (where we also use the equivalence 
(\ref{eqn:8.4}) from Remark 6.1). Hence, by Proposition 6.10, there exists 
$\pi \in NC(n)$ such that $\theta \leqleq \piodd \UU \kpeven$, and 
the latter inequality is clearly equivalent to having that
$\sigma \leqleq \pi$ and $\tau \leqleq K( \pi )$. Observe moreover 
that here we also obtain the uniqueness of $\pi$: if $\rho \in NC(n)$
is such that $\sigma \leqleq \rho$ and $\tau \leqleq K( \rho )$, then
it follows that $\theta \leqleq \rhodd \UU \kroeven$, and the uniqueness
part in Proposition 6.10 implies that $\rho = \pi$.  {\bf QED}

$\ $

$\ $

\setcounter{section}{7}
\begin{center}
{\large\bf 7. Proofs of Theorems 2 and 2'}
\end{center}

{\bf 7.1 Remark.}
The calculations made in this section will use the explicit
\setcounter{equation}{0}
formula for the coefficients of a boxed convolution $f \ \freestar \ g$,
where $f$ and $g$ are series in $\ncserk$. This formula says that
for every $n \geq 1$ and every $1 \leq i_1, \ldots , i_n \leq k$ we have
\begin{equation}  \label{eqn:9.1}
\cf_{(i_1, \ldots , i_n)}  (f \ \freestar \ g) =
\sum_{\pi \in NC(n)}
\cf_{(i_1, \ldots , i_n); \pi} (f) \cdot 
\cf_{(i_1, \ldots , i_n); K( \pi )} (g) .
\end{equation}
One can moreover extend Equation (\ref{eqn:9.1}) to a formula
which describes the generalized coefficients 
$\cf_{(i_1, \ldots , i_n); \pi}$ (as introduced in Notation 3.2.3) for 
the series $f \ \freestar \ g$. This says that
for every $n \geq 1$, every $1 \leq i_1, \ldots , i_n \leq k$, and 
every $\rho \in NC(n)$ we have 
\begin{equation}  \label{eqn:9.2}
\cf_{(i_1, \ldots , i_n); \rho}  (f \ \freestar \ g) =
\sum_{ \begin{array}{c}
{\scriptstyle \pi \in NC(n)}, \\
{\scriptstyle \pi \leq  \rho} 
\end{array} } \cf_{(i_1, \ldots , i_n); \pi} (f) \cdot 
\cf_{(i_1, \ldots , i_n); K_{\rho} ( \pi )} (g) .
\end{equation}
For a discussion of how one arrives to these formulas, we refer to
\cite{NS06}, Lecture 17.

$\ $

Instead of proving directly the Theorems 2 and 2', we will prove yet 
another equivalent reformulation of these theorems, which is stated 
as follows.

$\ $

{\bf 7.2 Theorem.} Let $f,g$ be two series in $\ncserk$. Then
\begin{equation}  \label{eqn:9.3}
\Reta ( f \ \freestar \ g ) = \Reta (f) \ \freestar \ \Reta (g).
\end{equation}

$\ $

{\bf Proof.} We fix $n \geq 1$ and $1 \leq i_1, \ldots , i_n \leq k$, and
we will prove the equality of the coefficients of $z_{i_1} \cdots z_{i_n}$
of the series appearing on the two sides of (\ref{eqn:9.3}). On the 
left-hand side we have:
\[
\cf_{(i_1, \ldots , i_n)} ( \Reta ( f \ \freestar \ g )) 
= \sum_{  \begin{array}{c}
{\scriptstyle  \rho \in NC(n),}  \\
{\scriptstyle  \rho \leqleq 1_n}  
\end{array} } \ \ \cf_{(i_1, \ldots , i_n); \rho} (f \ \freestar \ g)
\ \mbox{ (by Proposition 3.9.1)}
\]
\[
= \ \sum_{  \begin{array}{c}
{\scriptstyle  \rho \in NC(n),}  \\
{\scriptstyle  \rho \leqleq 1_n}  
\end{array} } \ \Bigl( \sum_{ \begin{array}{c}
{\scriptstyle  \sigma \in NC(n),}  \\
{\scriptstyle  \sigma \leq \rho}  
\end{array} } \ \ \cf_{(i_1, \ldots , i_n); \sigma} (f) \cdot
\cf_{(i_1, \ldots , i_n); K_{\rho} ( \sigma )} (g)  \Bigr)
\ \mbox{ (by Eqn.(\ref{eqn:9.2}))}
\]
\[
= \ \sum_{\sigma , \tau \in NC(n)} N'( \sigma , \tau)
\cf_{(i_1, \ldots ,i_n); \sigma} (f) \cdot 
\cf_{(i_1, \ldots ,i_n); \tau} (g),
\]
where for $\sigma, \tau \in NC(n)$ we denoted
\begin{equation}  \label{eqn:9.4}
N'( \sigma , \tau) := \card \{  \rho \in NC(n) \ | \ 
\rho \leqleq 1_n, \ \sigma \leq \rho, \ K_{\rho} ( \sigma ) = \tau \} .
\end{equation}
On the right-hand side of Equation (\ref{eqn:9.3}) the corresponding 
coefficient is:
\[
\cf_{(i_1, \ldots , i_n)} ( \Reta (f) \ \freestar \ \Reta (g) ) 
= \sum_{\pi \in NC(n)} \cf_{(i_1, \ldots , i_n); \pi} ( \Reta (f)) \cdot
\cf_{(i_1, \ldots , i_n); K( \pi )} ( \Reta (g)) 
\]
\[
= \ \sum_{\pi \in NC(n)} 
\Bigl( \sum_{ \begin{array}{c}
{\scriptstyle \sigma \in NC(n),}  \\
{\scriptstyle \sigma \leqleq \pi} 
\end{array} } \ \cf_{(i_1, \ldots , i_n); \sigma} (f) \Bigr) \cdot 
\Bigl( \sum_{ \begin{array}{c}
{\scriptstyle \tau \in NC(n),}  \\
{\scriptstyle \tau \leqleq K( \pi )} 
\end{array} } \ \cf_{(i_1, \ldots , i_n); \tau} (g) \Bigr) 
\ \mbox{ (by Eqn. (\ref{eqn:4.60}))}
\]
\[
= \ \sum_{\sigma , \tau \in NC(n)} N''( \sigma , \tau)
\cf_{(i_1, \ldots ,i_n); \sigma} (f) \cdot 
\cf_{(i_1, \ldots ,i_n); \tau} (g),
\]
where for $\sigma, \tau \in NC(n)$ we denoted
\begin{equation}  \label{eqn:9.5}
N''( \sigma , \tau) := \card \{  \pi \in NC(n) \ | \ 
\sigma \leqleq \pi, \ \tau \leqleq K( \pi ) \} .
\end{equation}
From the above calculations it is clear that (\ref{eqn:9.3}) will follow
if we can prove that
\[
N' (\sigma, \tau ) = 
N'' (\sigma, \tau ), \ \ \forall \, \sigma , \tau \in NC(n).
\]
Now, the content of Corollary 6.11 is that
\begin{equation}  \label{eqn:9.6}
N''( \sigma , \tau ) = \left\{  \begin{array}{ll}
1, & \mbox{if $\tau \leq K( \sigma )$ and 
           $P_{\tau}(n) = P_{\sigma}^{-1} (1)$ }     \\ 
0, & \mbox{otherwise.}
\end{array}  \right.
\end{equation}
So it remains to prove that $N' (\sigma , \tau )$ is also described 
by the right-hand side of (\ref{eqn:9.6}).

Let us observe that always $N' (\sigma , \tau ) \in \{ 0,1 \}$. Indeed,
if there exists $\rho \in NC(n)$ such that $\rho \geq \sigma$ and 
$K_{\rho} (\sigma) = \tau$,
then $\rho$ is uniquely determined -- this is because the associated 
permutation $P_{\rho}$ is determined, $P_{\rho} = P_{\sigma} P_{\tau}$.
In order to prove that $N' ( \sigma , \tau )$ is equal to the quantity
on the right-hand side of (\ref{eqn:9.6}) it will therefore suffice to
verify the following equivalence:
\begin{equation}   \label{eqn:9.7}
\Bigl( \ \tau \leq K( \sigma ) \mbox{ and }
           P_{\tau}(n) = P_{\sigma}^{-1} (1) \ \Bigr) \ \Leftrightarrow \
\Bigl( \ \begin{array}{c}
\exists \, \rho \in NC(n) \mbox{ such that}    \\
\rho \leqleq 1_n, \ \sigma \leq \rho, \mbox{ and } K_{\rho} (\sigma) = \tau
\end{array} \ \Bigr).
\end{equation}

\vspace{6pt}

{\em Verification of ``$\Rightarrow$'' in (\ref{eqn:9.7}).} Consider the 
relative Kreweras complement of $\tau$ in $K( \sigma )$, and then 
consider the partition 
\[
\rho := K^{-1} \bigl(  K_{ K( \sigma ) } ( \tau ) \bigr) \in NC(n).
\]
From Equation (\ref{eqn:2.3}) in Remark 2.3 we have that 
$K_{K( \sigma )} ( \tau ) \leq K( \tau )$; if we then apply the 
order-reversing map $K^{-1}$ to both sides of this inequality, we get
that $\rho \geq \sigma$. An 
immediate calculation involving the permutations associated to 
$\rho, \sigma$ and $\tau$ gives us that $P_{\rho} = P_{\sigma} P_{\tau}$,
and this in turn implies that $K_{\rho} ( \sigma ) = \tau$. Finally,
observe that $P_{\rho} (n) = P_{\sigma} P_{\tau} (n) = 1$ (with the 
latter inequality following from the fact that 
$P_{\sigma}^{-1} ( 1) = P_{\tau} (n)$). This shows that
$\rho \leqleq 1_n$, and completes this verification.

\vspace{6pt}

{\em Verification of ``$\Leftarrow$'' in (\ref{eqn:9.7}).} 
Let $\rho \in NC (n)$ be such that $\rho \leqleq 1_n$, 
$\rho \geq \sigma$, and $K_{\rho} (\sigma) = \tau$. From Equation 
(\ref{eqn:2.3}) in Remark 2.3 we obtain that
$\tau = K_{\rho} ( \sigma ) \leq K( \sigma )$. On the other hand
the permutations associated to $\rho, \sigma, \tau$ satisfy
$P_{\tau} = P_{\sigma}^{-1} P_{\rho}$  (because 
$\tau = K_{\rho} ( \sigma )$), and $P_{\rho} (n)=1$ (because 
$\rho \leqleq 1_n$). Hence we have 
$P_{\tau} (n) = P_{\sigma}^{-1} ( P_{\rho} (n)) =  
P_{\sigma}^{-1} (1)$, as required.  {\bf QED}

$\ $

{\bf 7.3 Remark} {\em (proofs of Theorems 2 and 2').} In the introduction 
section it was shown how Theorem 2 is derived from Theorem 2', and here 
we show how Theorem 2' follows from the above Theorem 7.2. Let $\mu$ and 
$\nu$ be two distributions from $\Dalg (k)$. Consider the formula 
(\ref{eqn:1.13}) which is satisfied by $\mu$ and $\nu$, and apply 
$\Reta$ to both its sides. We obtain
\begin{align*}
\Reta \Bigr( R_{\mu \boxtimes \nu} \Bigl) 
& = \Reta \Bigr( R_{\mu} \ \freestar \ R_{\nu} \Bigl)         \\
& = \ \Reta ( R_{\mu} ) \ \freestar \ \Reta ( R_{\nu} )  \ \
\mbox{ (by Theorem 7.2).}
\end{align*}
Since $\Reta$ maps an $R$-transform to the $\eta$-series of the same 
distribution, we have thus obtained that $\eta_{\mu \boxtimes \nu}$ =
$\eta_{\mu} \ \freestar \ \eta_{\nu}$, as stated in Theorem 2'.

$\ $

$\ $

$\ $ 

$\ $

Serban T. Belinschi and Alexandru Nica

Department of Pure Mathematics,
University of Waterloo

Waterloo, Ontario N2L 3G1, Canada

Email: sbelinsc@math.uwaterloo.ca, anica@math.uwaterloo.ca

\end{document}